\documentclass{article}
\usepackage[utf8]{inputenc}
\usepackage[cmintegrals]{newtxmath}
\usepackage[T1]{fontenc}
\usepackage[top=0.75in,bottom=0.75in,left=1.0in,right=1.0in]{geometry}

\usepackage{graphicx}
\usepackage{mathrsfs}
\usepackage{amsmath,stackrel}
\usepackage{graphicx}
\usepackage{multirow}
\usepackage{mathtools}
\usepackage{ marvosym }
\usepackage{tikz-cd}
\usepackage{pst-node}

\newcommand{\so}{\mathfrak{so}(3)}

\newcommand{\R}{\mathbb{R}}

\newcommand{\SO}{\text{SO}}
\newcommand{\eff}{\text{eff}}

\usepackage{bm}
\usepackage[colorlinks=true]{hyperref}
\hypersetup{urlcolor=blue, citecolor=red}
\newtheorem{theorem}{Theorem}[section]
\newtheorem{corollary}{Corollary}
\newtheorem{lemma}[theorem]{Lemma}
\newtheorem{proposition}{Proposition}

\newtheorem{definition}[theorem]{Definition}
\newtheorem{remark}{Remark}

\title{Geometric Control of two Quadrotors Carrying a Rigid Rod with Elastic Cables}
\author{Jacob R. Goodman and Leonardo J. Colombo}
\date{}

\begin{document}

\maketitle

\begin{abstract}
This paper presents the design of a geometric trajectory tracking controller for the cooperative task of two quadrotor UAVs (unmanned aerial vehicles) carrying and transporting a rigid bar, which is attached to the quadrotors via inflexible elastic cables. The elasticity of the cables together with techniques of singular perturbation allows a reduction in the model to that of a similar model with inelastic cables. In this reduced model, we design a controller such that the rod exponentially tracks a given desired trajectory for its position and attitude, under some assumptions on initial error. 
We then show that exponential tracking in the reduced model  corresponds to exponential tracking of the original elastic model. We also show that the previously defined control scheme provides uniform ultimate boundedness in the presence of unstructured bounded disturbances.
\end{abstract}


\section{Introduction}

The use of aerial robots has become increasingly popular in the last decades due to their superior mobility and versatility in individual and cooperative tasks \cite{PalCruIRAM12}. For instance, aerial robots equipped with manipulators can be utilized for mobile manipulation tasks such as rescue operations and transportation. Recently, control design of multiple aerial robots transporting objects has been studied in the literature \cite{MazKonJIRS10}, \cite{MicFinAR11}, \cite{pedro1}, \cite{pedro2}, \cite{jose}.  Several of these aerial robots can be used to transport heavier payloads thus expanding the capabilities of a single aerial robot \cite{gas}. 

 In aerial transportation, a cable establishes a physical connection between the UAV and the cargo. Geometric nonlinear controllers of multiple quadrotors with a suspended point-mass load were studied in \cite{SreLeePICDC13} and with a rigid body load in \cite{Lee-CDC}, \cite{wu}. In \cite{god}, the authors model the cables as flexible chains comprised of inflexible links with mass. While these works have considered the cable to be inelastic, we are instead motivated by applications where the elasticity in the cable tethers cannot be ignored without compromising the validity of the estimations \cite{jo}. A single quadrotor carrying a point mass payload with elastic cables has been studied in \cite{Elastic} and \cite{wu}. In this paper we study the problem of two quadrotors transporting a rigid rod suspended through elastic cables. To the best of our knowledge, it is the first time that a geometric controller for such a transportation task with quadrotors has been developed. 

The use of geometric controllers in the UAVs literature has been extensively developed in the last years (see, for instance,  \cite{amit1}, \cite{amit2}, \cite{amit3}, \cite{amit4}, \cite{Lee-CDC}, \cite{Lee-TAC}, \cite{Lee-TCST}, \cite{LeeSrePICDC13}, and references therein). In this paper, we propose a coordinate-free form of the equations of
motion for the cooperative task, which are derived according to Lagrangian mechanics on manifolds—in particular, via the Lagrange d-Alembert principle for forced systems \cite{abloch}. Working directly on the manifold allows us to avoid potential
singularities of local parameterizations (e.g. Euler Angles), generating agile maneuvers of the payload in a uniform manner. In particular, in this work, a geometric control scheme taking the form of a feedback linearization together with a geometric PD controller, as in  \cite{muga} and \cite{muga2}, is designed such that the rigid bar exponentially reach and follows a given desired trajectory of both the bar's position and attitude. 

The main contributions of this work are: (i) the modelling and subsequent derivation of the corresponding equations of motion for the cooperative task of two quadrotor UAVs transporting a rigid rod via inflexible elastic cables. The modelling and dynamics are summarized in Proposition \ref{equationsprop}. (ii) Reduction of these equations of motion to the case of inelastic cables under the assumption of sufficiently high spring damping and stiffness. This can be considered an extension of the results obtained in \cite{Elastic}, which studies the case of a single quadrotor transporting a point mass load via an inflexible elastic cable. This is developed in Section \ref{sec4} and employed in Section \ref{sec5} where (iii) we provide a geometric control scheme for the exponential tracking of the load position and attitude to some desired trajectories. The proposed controller is inspired by, and takes a similar form to, that which is found in \cite{Lee-TCST}—which studies the case of a team of quadrotor UAVs transporting a rigid body with rigid cables. It is important to note that such a controller cannot be applied directly to our problem, as it requires a minimum of $3$ UAVs and at least a $2$-dimensional rigid body. The major consequences of the rigid rod being a $1$-dimensional rigid body come in the description of the attitude - which lies in the unit sphere $S^2$ as opposed to the special orthogonal group $\SO(3)$ — and also by the fact that the inertia tensor is \textit{singular}. Lyapunov analysis is used to determine sufficient conditions for exponential tracking. Theorem \ref{gainsth} then proves the existence of stabilizing gains—for sufficiently small initial errors in the cables—which satisfy the conditions. Such a proof was not previously seen in the literature, and in principle yields some insight into the relationship between gains. Finally, we handle the case of unstructured bounded disturbances acting on our system, which also had not been seen in the literature. In particular, (iv) in Theorem \ref{gainsthdist} we show that the same control scheme will yield uniform ultimate bounds in the case of unstructured bounded disturbances acting on the system. Moreover, the ultimate bound can be made arbitrarily small by choosing gains appropriately.

The rest of the paper is structured as follows. In Section \ref{sec3} we model and derive the dynamical system describing the task of carrying and transport a rigid rod between two quadrotors by elastic cables. This is done by constructing the Lagrangian of the cooperative system and subsequently applying the Lagrange d'Alemebrt principle. Section \ref{sec4} reduces the dynamical model introduced in Section \ref{sec3} by employing singular perturbation theory techniques. The main results of the work are given in Section \ref{sec5} and \ref{sec_dist}. We first introduce configuration error functions for each state variable, from which we derive the error dynamics. The geometric controls are constructed in order that the origin is an exponentially stable equilibrium point of the error dynamics. After this analysis we return to the original model to adequate our geometric control such that trajectories of the original model exponentially track their desired trajectories up to a neighborhood whose size shrinks uniformly with increasing spring stiffness and damping. In section \ref{sec_dist}, we introduce unstructured bounded disturbances to the reduced model and show that the same control scheme can be applied to achieve uniform ultimate boundedness. Numerical simulations are shown to validate the theoretical results. 



\begin{section}{Modelling and Control Equations}\label{sec3}

In this section we model and derive the dynamical system describing the cooperative task between the quadrotors. This can be done by constructing the total kinetic and potential energies of the mechanical system describing the cooperative task—in addition to the virtual work done by non-conservative forces—and subsequently using the tools of Lagrangian Mechanics on manifolds \cite{abloch}. 

Consider two identical quadrotor UAVs transporting a rigid rod of length $2L_r$ and total mass $m_r$. The rod is considered inflexible and of uniform mass density. The endpoints of the rod are connected to the center of mass of each quadrotor via a massless inflexible elastic cable of rest length $L_c$, as it is shown in Figure \ref{figuav}.

\begin{figure}[h!]
\begin{center}
 \includegraphics[width=8.cm]{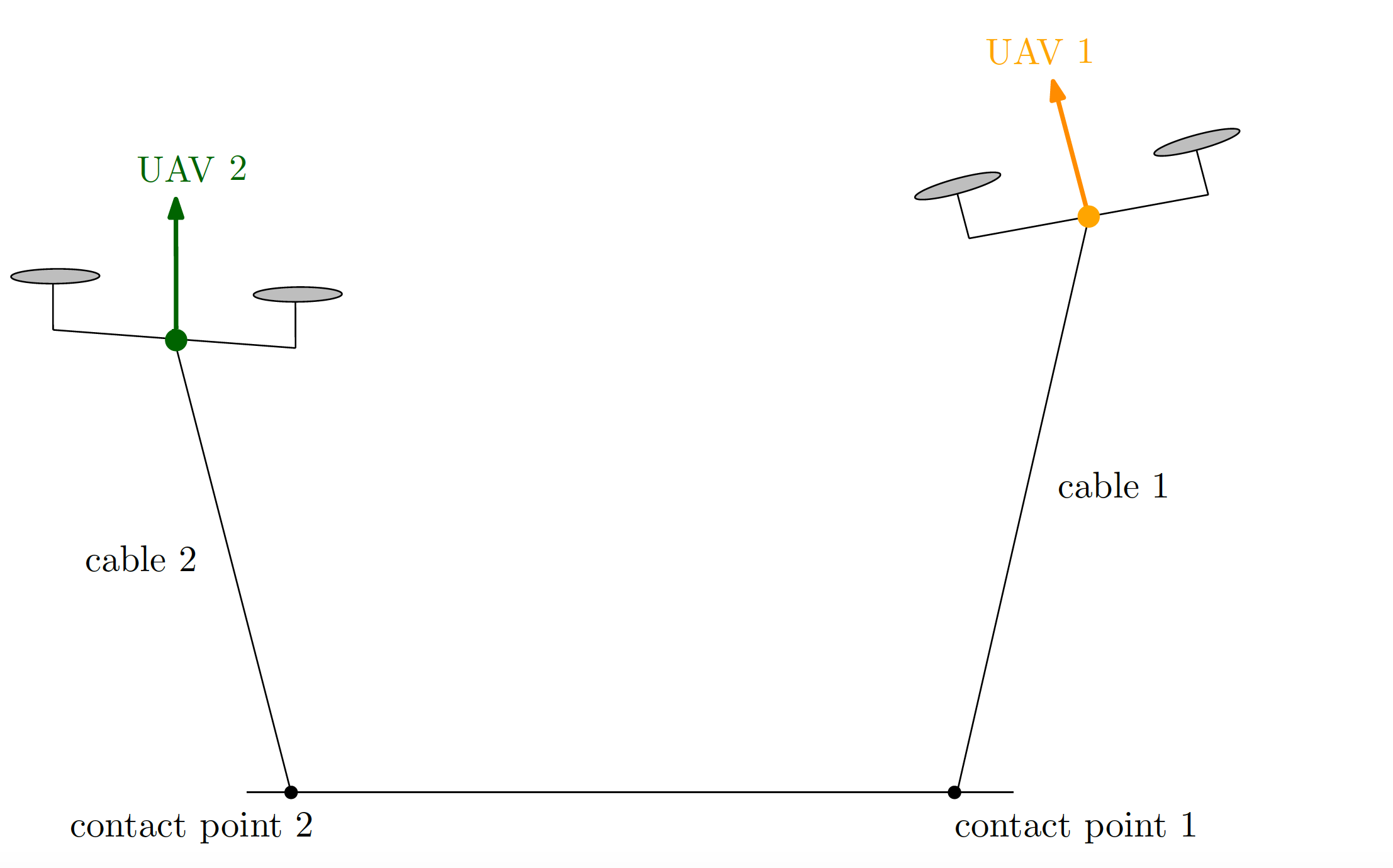}
 \caption{Modeling of the mechanical system describing the cooperative transportation task.}
 \label{figuav}
 \end{center}
\end{figure}


The basic notation and methodology is fairly
standard within the literature and we have attempted to use traditional symbols and definitions wherever feasible. Table \ref{table} provides the symbols and geometric spaces that are used frequently throughout the paper.

\begin{table}[h!]
\begin{center}
\begin{tabular}{|l|c|p{0.7\linewidth}|}
\textbf{Variable} & \textbf{Space} & \textbf{Description} \\
\hline
$m_Q$ & $\mathbb{R}$ & Mass of quadrotor. \\
$m_r$ & $\mathbb{R}$ & Mass of rod. \\
$2L_r$ & $\mathbb{R}$ & Length of rod. \\
$L_c$ & $\mathbb{R}$ & Rest length elastic cables. \\
$x_{Q_j}$ & $\mathbb{R}^3$ & Position of quadrotor $j$ in inertial frame. \\
$q_r$ & $S^2$ & Attitude of rod in inertial frame. \\
$x_r$ & $\mathbb{R}^3$ & Position of center of mass of rod in inertial frame. \\
$q_j$ & $S^2$ & Position vector of cable suspended from quadrotor $j$. \\
$l_j$ & $\mathbb{R}$ & Length of elastic cable attached to quadrotor $j$. \\
$R_j$ & $\text{SO}(3)$ & Attitude of quadrotor $j$. \\
$J_Q$ & $\hbox{Sym}_{\succ 0}(\R)$ & Moment of inertia of quadrotor.\\
$J_r$ & $\hbox{Sym}_{\succeq 0}( \R)$ & Moment of inertia of the rod.\\
$\omega_r$ & $TS^2$ & Angular velocity of rod in inertial frame. \\
$v_r$ & $\mathbb{R}^3$ & Translational velocity of center of mass of rod in inertial frame. \\
$\omega_j$ & $TS^2$ & Angular velocity of cable $j$ in the inertial frame. \\
$\Omega_j$ & $\so(3)$ & Angular velocity of quadrotor $j$ in body frame. \\
$u_j$ & $\mathbb{R}^3$ & Net thrust applied vertically in the body frame of quadrotor $j$. \\
$M_j$ & $\mathbb{R}^3$ & Moment vector in the body frame of quadrotor $j$. \\
$u^{\perp_j}$ & $\mathbb{R}^3$ & The component of $u$ that is perpendicular to $q_j$. \\
$u^{\parallel_j}$ & $\mathbb{R}^3$ & The component of $u$ that is parallel to $q_j$. \\

\end{tabular}
\caption{Nomenclature}
\label{table}
\end{center}
\end{table}

The configuration space of the mechanical system is given by \\$Q = \underbrace{(S^2 \times \mathbb{R}^3)}_\text{Rigid rod} \times \underbrace{(S^2 \times \R)^2}_\text{Cables} \times \underbrace{(\SO(3) \times \SO(3))}_\text{Quadrotor attitudes}$. We fix an inertial frame such that $e_3$ is oriented opposite to the direction of gravitational acceleration, and define the following auxiliary variables to describe the position of the quadrotors in this inertial frame $x_{Q_j} = x_r + (-1)^j L_r q_r - l_j q_j$ for $j = 1,2$. With these coordinates, the translational kinetic energy of each quadrotor can be described by $\frac12 m_Q ||\dot{x}_{Q_j}||^2$, where $m_Q$ denotes the mass of the quadrotor. Similarly, the translational kinetic energy of the rod is $\frac{1}{2}m_r ||\dot{x}_r||^2$. As the quadrotors and rod are rigid bodies, we further have rotational kinetic energy components in the total kinetic energy. Fixing a body frame to each quadrotor and denoting the angular velocity in this body frame by $\Omega_j\in\mathbb{R}^{3}$, the angular kinetic energy is given by $\frac{1}{2}\Omega_j^T J_Q \Omega_j$, where $J_Q$ is a symmetric positive-definite inertia tensor. $\Omega_j$ is defined implicitly by the kinematic equation $\dot{R}_j = R_j \hat{\Omega}_j$, where $\hat{\cdot}: \mathbb{R}^3 \to \so(3)$ is the \textit{hat isomorphism} which maps vectors on $\mathbb{R}^{3}$ to $(3\times 3)$ skew-symmetric matrices

\begin{align*}
\Omega = \begin{bmatrix} \Omega_1 \\ \Omega_2 \\ \Omega_3 \end{bmatrix} \mapsto \begin{bmatrix} 0 & -\Omega_3 & \Omega_2 \\ \Omega_3 & 0 & -\Omega_1 \\ -\Omega_2 & \Omega_1 & 0 \end{bmatrix} := \hat{\Omega}.
\end{align*}

The angular kinetic energy of the rod is similarly given by $\frac12 \Omega_r^T J_r \Omega_r$, however the symmetry of the rod allow us to reinterpret the expression. In particular, we will fix a body frame to one of the endpoints of the rod with $e_1 := q_r$ as one of the (orthonormal) basis vectors. In this frame, the component of the (body) angular velocity along $e_1$ is identically zero. Moreover, we can describe the inertia tensor explicitly by $J_r= \begin{bmatrix} 0 & 0 & 0 \\ 0 & \frac43 m_r L_r^2 & 0 \\ 0 & 0 & \frac43 m_r L_r^2 \end{bmatrix}$. From this is clear that $\frac12 \Omega_r^T J_r \Omega_r = \frac23 m_r L_r^2 ||\Omega_r||^2$. 

The angular velocity of the rod in the inertial frame is defined implicitly by the kinematic relation $\dot{q}_r = \omega_r \times q_r$ together with the condition that $\omega_r^T q_r = 0$, from which it is obvious that that $||\dot{q}_r|| = ||\omega_r||$. Moreover, for some rotation $R_r \in \SO(3)$, we have $\Omega_r = R_r \omega_r$, and since the norm of a vector is invariant under the action of $\SO(3)$, we see that $||\Omega_r|| = ||\dot{q}_r||$. Altogether, the angular kinetic energy of the rod is given by $\frac23 m_r L_r^2 ||\dot{q}_r||^2$.

Recall that elements in the tangent space $T_{R}SO(3)$ are identified with elements in $SO(3)\times\mathfrak{so}(3)$ by a left-trivialization. That is, the diffeomorphism $(R,\dot{R})\in T_{R}SO(3)\mapsto (R,R^{-1}\dot{R})=:(R,\hat{\Omega})\in SO(3)\times\mathfrak{so}(3)$ (see \cite{HSS} for details). Therefore, after a left trivialization of $TSO(3)$, the tangent bundle of $Q$ can be identified as
$$TQ \cong  (S^2 \times TS^2 \times \mathbb{R}^3\times \mathbb{R}^3)\times(S^2 \times TS^2 \times \mathbb{R}\times \mathbb{R}) \times  (\SO(3)\times \so(3))^2.$$

Finally, the total kinetic energy $K: TQ \to \R$ of the system is given by summing the respective translational and angular kinetic energies of the quadrotors and rigid rod:

$$K =\underbrace{\frac{1}{2}m_r ||\dot{x}_r||^2 + \sum_{j=1}^{2} \frac{1}{2}m_Q ||\dot{x}_{Q_j}||^2}_\text{Translational K.E.} + \underbrace{\frac23 m_r L_r^2 ||\dot{q}_r||^2 + \sum_{j=1}^{2} \frac{1}{2}\Omega_j^T J_Q \Omega_j. }_\text{Angular K.E.}$$ 

Moreover, the total potential energy $U: Q \to \R$ of the system is given by

$$U = \sum_{j=1}^{2}\underbrace{m_Q ge_3^T x_{Q_j} + m_r g e_3^T x_r}_\text{Gravitational P.E.} + \sum_{j=1}^{2} \underbrace{ \frac{1}{2}k(L_c - l_j)^2,}_\text{Elastic P.E.}$$which corresponds to the gravitational potential energies of the quadrotors and the rigid rod, as well as the elastic potential of the cables. As usual, the Lagrangian of the system $L: TQ \to \R$ is defined by $L:= K - U$.

Control inputs for each quadrotor are denoted by $u_j, M_j \in \R^3$. The first is a thrust controller corresponding to the total lift force exerted on the quadrotor by the spinning propellers. In particular, $u_j=f_jR_je_3\in\mathbb{R}^3$, where $f_j\in\mathbb{R}$ is the total thrust magnitude and $e_3 =[0,\, 0,\,1]^{T}\in\mathbb{R}^3$. The second is a moment controller, which is related to the torque induced on the quadrotor by propellers. Alternatively, one may choose to control the total thrust of \textit{each} propeller individually. However, we opt for the former approach because it leads nicely to the separation of the quadrotor's attitude dynamics from the rest of the system's dynamics. The thrust generated by the $i$-th propeller along
the $e_3$ axis can be determined by the total thrust and the moment controller as in \cite{LLM}.

Note that these controls take the form of non-conservative external forces, so that we must use the Lagrange d'Alembert Variational Principle (see \cite{abloch} for instance) -- with controls playing the role of the virtual forces in our system -- to obtain our system dynamics from the Lagrangian $L$. We further wish to add a non-conservative force corresponding to a damping in the elastic cables. That is, a velocity dependent force that serves to reduce the amplitude of oscillations in our elastic cable. In particular, we will opt to make this force \textit{proportional} to the velocity, with constant of proportionality $c > 0$. 

Denote by $C^{\infty}(Q,q_0,q_T)$ the space of smooth function from $[0,T]$ to $Q$ with fixed endpoints points, denoted by $q_0$ and $q_T$, respectively. Consider the action functional $\mathcal{A}:C^{\infty}(Q,q_0,q_T)\to\mathbb{R}$ given by \begin{equation}\label{action}
\mathcal{A}(c(t))=\int_{0}^{T} L(c(t),\dot{c}(t))\, dt + \sum_{j=1}^{2} \int_{0}^{T} \left(||f_jR_je_3||_{\mathbb{R}^{3}}^{2} + ||\hat{M}_j||^2_{\mathfrak{so}(3)} - c \dot{l}_j\right)\,\, dt = 0,
\end{equation}  where $||\hat{M}_j||_{\mathfrak{so}(3)}:=\langle\hat{M}_j,\hat{M}_j\rangle^{1/2}=\sqrt{\hbox{Tr}(\hat{M}_j^{T}\hat{M}_j)}$, with  $c(t):=(q_r(t),x_r(t),q_j(t),l_j(t),R_j(t)))\in C^{\infty}(Q,q_0,q_T)$.

In order to use the aforementioned Lagrange d'Alembert Variational Principle, we must describe the variations of our state variables. These variations must be tangent vectors in the tangent spaces of the submanifolds of the configuration space in which the state variables live. In addition, they must vanish at the end points, because tangent vectors on the tangent bundle of $C^{\infty}(Q,q_0,q_T)$ must satisfy such a condition (see for instance \cite{HSS}, \cite{MR}). 

In particular, we choose
$\delta x_r \in \mathbb{R}^3$ and $\delta l_j \in \mathbb{R}$  arbitrary,  $\delta q_j = \frac{d}{d\epsilon}\mid_{\epsilon=0}\exp(\epsilon \hat{\xi}_j)q_j=\xi_j \times q_j \in T_{q_j} S^2$, $\delta q_r=\frac{d}{d\epsilon}\mid_{\epsilon=0}\exp(\epsilon \hat{\xi}_r)q_r=  \xi_r \times q_r\in T_{q_r} S^2$, satisfying $\xi_j\cdot q_j=0$ and $\xi_r\cdot q_r=0$, for arbitrary vectors $\xi_j, \ \xi_r \in \mathbb{R}^3$ and $j = 1,2$. In addition by defining the curve on the Lie algebra $\mathfrak{so}(3)$ given by $\hat{\eta}_j = R_j^T \delta R_j \in \so(3)$, it can be shown that (see for instance \cite{MR} Chapter 13) $\widehat{\delta \Omega}_j = \hat{\Omega}_j \hat{\eta}_j + \dot{\hat{\eta}}_j$ with $\hat{\eta}_j$ satisfying $ \hat{\eta}_j(0)= \hat{\eta}_j(T)=0$ (since $\delta R_j(0)=\delta R_j(T)=0)$. Moreover, we have the following relations
\begin{align*}
\delta x_{Q_j} &= \delta x_r + (-1)^jL_r \delta q_r - (\delta l_j) q_j - l_j (\delta q_j), \\
\delta \dot{x}_{Q_j} &= \delta \dot{x}_r + (-1)^j L_r \delta \dot{q}_r - (\delta \dot{l}_j) q_j - (\delta l_j) \dot{q}_j - \dot{l}_j (\delta q_j) - l_j (\delta \dot{q}_j). 
\end{align*}

\begin{proposition}\label{equationsprop}
Critical points of the action functional $\mathcal{A}$ for variations with fixed endpoints corresponds with solutions of the controlled Euler-Lagrange equations
\begin{align}
\dot{x}_r &= v_r, \label{xkin} \\
m_{\eff} ( \dot{v}_r + ge_3) &= u_1 + u_2 + m_Q (\ddot{\zeta}_1 + \ddot{\zeta}_2 ), \label{xdyn}\\
\dot{q}_r &= \omega_r \times q_r, \label{qrkin} \\
 I_{\eff} \dot{\omega}_r &= q_r \times \left[ u_2 - u_1 + m_Q( \ddot{\zeta}_{2} - \ddot{\zeta}_1) \right], \label{qrdyn} \\
m_Q q_j^T \ddot{\zeta}_j &= m_Q q_j^T ( \dot{v}_r + (-1)^j L_r \ddot{q}_r + ge_3 - \frac1{m_Q}u_j) - c \dot{l}_j  + k(L_c - l_j), \label{ldyn}\\
q_j \times \ddot{\zeta}_j &= q_j \times ( \dot{v}_r + (-1)^j L_r \ddot{q}_r + ge_3 - \frac1{m_Q} u_j),\label{qjdyn} \\
\dot{R}_j &= R_j \hat{\Omega}_j, \label{Rkin}\\
J_Q \dot{\Omega}_j &= J_Q \Omega_j \times \Omega_j + M_j,  \quad \text{for $j = 1,2$,} \label{Rdyn}
\end{align} where $m_{\eff} := 2m_Q + m_r, \ I_{\eff} := (2m_Q + \frac23 m_r)L_r$, and $\zeta_{j} := l_j q_j$.
\end{proposition}
\begin{remark}
Equations \eqref{xkin}-\eqref{xdyn} describe the kinematics and dynamics of the rod's position, respectively. Similarly, equations \eqref{qrkin}-\eqref{qrdyn} describe the kinematics and dynamics of the rod's attitude, and equations \eqref{Rkin}-\eqref{Rdyn}, corresponds with the kinematics and dynamics of the attitudes of each quadrotor, respectively. 

Equation \eqref{ldyn}, indexed for $j = 1,2$, describes the dynamics of the lengths of the elastic cables. This can be understood by observing that the projection of $\ddot{\zeta}_j$ onto $q_j$ preserves the acceleration of the length (which is inherently oriented along the cable), while removing the acceleration of the attitude from consideration with the identity $q_j^T \ddot{q}_j = -||\dot{q}_j||^2$. Conversely, equation \eqref{qjdyn}, indexed for $j = 1,2$, describes the dynamics of the attitudes of the elastic cables, as the cross-product with $q_j$ preserves the acceleration of the cable attitude while annihilating the acceleration of the cable length. 
\end{remark}

\textit{Proof of Proposition 2.1:} We wish to apply Lagrange- d'Alembert Variational Principle. Therefore, our system dynamics must satisfy

\begin{equation}\label{actionvariation}
\delta \int_{0}^{T} L(c(t),\dot{c}(t))\, dt + \sum_{j=1}^{2} \int_{0}^{T} \left(\delta x_{Q_j}^T u_j + \left<R_j^T \delta R_j, \hat{M}_j\right> - c \dot{l}_j\delta l_j\right)dt = 0,
\end{equation} where the integral on the right represents the virtual work done by the thrust controls $u_j$, the moment controls $M_j \in \R^3$, and the spring damping, respectively. 

Expanding the variations within \eqref{actionvariation}, substituting the corresponding infinitesimal variations, and grouping like terms, we obtain
\begin{align*}
0=& \int_0^{T} \left[ \delta \dot{x}_r^T ( m_{\eff} \dot{x}_r - m_Q(\dot{\zeta}_1 + \dot{\zeta}_2)) + \delta x_r^T (-m_{\eff} ge_3 + u_1 + u_2 ) \right] dt \\
+& \  L_r  \int_0^{T} \left[ \dot{\xi}_r^T ( q_r \times (I_{\eff} \dot{q}_r + m_Q(\dot{\zeta}_2 - \dot{\zeta}_1 )) + \xi_r^T (\dot{q}_r \times m_Q(\dot{\zeta}_2 - \dot{\zeta}_1) + q_r \times (u_2 - u_1) \right]dt \\
-& \sum_{j=1}^{2}\int_0^{T} \left[ m_Q(\delta \dot{l}_j) q_j^T \dot{x}_{Q_j} + \delta l_j ( m_Q \dot{q}_j^T \dot{x}_{Q_j} - m_Q ge_3^T q_j  - k(L_c - l_j) + c \dot{l}_j + q_j^T u_j ) \right] dt \\
-& \sum_{j=1}^{2}\int_0^{T} \left[\xi_j^T \left( q_j \times (m_Q \dot{l}_j \dot{x}_{Q_j} - m_Q g l_j e_3 + l_j u_j) + \dot{q}_j \times m_Q l_j \dot{x}_{Q_j} \right) + \dot{\xi}_j^T (q_j \times m_Q l_j \dot{x}_{Q_j}) \right]dt \\
+&\sum_{j=1}^{2} \int_0^{T}  \left[ \dot{\eta}_j^T J_Q \Omega_j  + \eta_j^T \left(  J_Q \Omega_j \times \Omega_j + M_j \right) \right]dt,
\end{align*}where $\zeta_j := l_j q_j$,  $m_{\eff} := 2m_Q + m_r$, and $ I_{\eff} := (2m_Q + \frac23 m_r)L_r$. Integrating by parts and applying the equality of mixed partial derivatives and the fact that variations vanish on the endpoints, we obtain
 \begin{align*}
0=& \int_{0}^{T} \delta x_r \left[ m_{\eff} (\ddot{x}_r + ge_3) - m_Q(\ddot{\zeta}_1 + \ddot{\zeta}_2) - (u_1 + u_2 )) \right] dt + \sum_{j=1}^{2}\int_0^{T} \eta_j^T \left[ J_Q \dot{\Omega}_j - J_Q \Omega_j \times \Omega_j - M_j \right]dt \\
&+ \  L_r  \int_0^T \xi_r^T \left[  I_{\eff}(q_r \times \ddot{q}_r) - q_r \times \left( u_2 - u_1 + m_Q( \ddot{\zeta}_{2} - \ddot{\zeta}_1) \right) \right]dt \\
&- \sum_{j=1}^{2}\int_0^{T} \delta l_j \left[ m_Q q_j^T( \ddot{x}_r + (-1)^j L_r \ddot{q}_r - \ddot{\zeta}_j+ ge_3) - c \dot{l}_j  + k(L_c - l_j) - q_j^T u_j \right] dt \\
&- l_j \sum_{j=1}^{2}\int_0^{T} \xi_j^T \left[m_Q q_j \times ( \ddot{x}_r + (-1)^j L_r \ddot{q}_r - \ddot{\zeta}_j + ge_3) - q_j \times u_j \right]dt. \\
\end{align*}

Each of these integrals can be treated independently, as their respective variations are independent. That is, for the above equation to be satisfied, we necessarily have that each integral vanish identically. Applying the Fundamental Lemma of the Calculus of Variations \cite{GF} to each integral yields the dynamical system:
\begin{align*}
m_{\eff} ( \ddot{x}_r + ge_3) =& u_1 + u_2 + m_Q (\ddot{\zeta}_1 + \ddot{\zeta}_2 ),\\ I_{\eff}(q_r \times \ddot{q}_r) =& q_r \times \left[ u_2 - u_1 + m_Q( \ddot{\zeta}_{2} - \ddot{\zeta}_1) \right],\\
m_Q q_j^T( \ddot{x}_r + (-1)^j L_r \ddot{q}_r - \ddot{\zeta}_j+ ge_3) =& c \dot{l}_j  - k(L_c - l_j) + q_j^T u_j, \\
m_Q q_j \times ( \ddot{x}_r + (-1)^j L_r \ddot{q}_r - \ddot{\zeta}_j + ge_3) = &q_j \times u_j, \\
\dot{R}_j =&R_j \hat{\Omega}_j,\\ d J_Q \dot{\Omega}_j =& J_Q \Omega_j \times \Omega_j + M_j,  \qquad\text{for $j = 1,2$}.
\end{align*}
where we have made the assumption that $l_j \ne 0$. After implicitly defining the translational and angular velocities of the load with the kinematic equations $\dot{v}_r = x_r$ and $\dot{q}_r = \omega_r \times q_r$, and rearranging terms, we have the desired dynamical control system. \hfill$\square$


\end{section}


\begin{section}{Reduced Model}\label{sec4}

While the use of elastic cables provide the benefit of reducing impulsive forces on the bar, large or rapid oscillations of the bar can produce undesired aggressive movements, compromising the safety of the cooperative task. Therefore, the strategy is to use elastic cables with high stiffness and damping to guarantee the safety for the bar in the transportation task. To this end, we will employ techniques from singular perturbation theory \cite{SPT} to study such a situation. 


In particular, we will consider the case that $\displaystyle{k = \frac{\bar{k}}{\epsilon^2}}$ and $\displaystyle{c = \frac{\bar{c}}{\epsilon}}$ with $\bar{k}, \bar{c} > 0$ and $\epsilon > 0$ sufficiently small, and we will show that the dynamics approach that of the same model with inelastic cables (that is, with $l \equiv L_c$) as $\epsilon \to 0$. We further consider a change of variables of the form $l_j = \epsilon^2 y_j + L_c$ and $\dot{l}_j = \epsilon z_j$, which is motivated by observing that $k( L_c - l_j) = -\bar{k}y_j$ and $c \dot{l}_j = \bar{c} z_j$. From this, we can see that $\zeta_j = (\epsilon^2 y_j + L_c)q_j$. Therefore, $\ddot{\zeta}_j = L_c\ddot{q}_j + \epsilon (\dot{z}_j q_j + z_j \dot{q}_j) + \epsilon^2 y_j \ddot{q}_j$. Making these substitutions into the dynamics described in Proposition \ref{equationsprop}, in addition to defining the angular velocity of the cables $\omega_j$ by $\dot{q}_j = \omega_j \times q_j$ and $\omega_j^T q_j = 0$, we obtain
\begin{align*}
\dot{x}_r =& v_r, \\
m_{\eff} ( \dot{v}_r + ge_3) =& u_1 + u_2 + m_Q L_c (\ddot{q}_1 + \ddot{q}_2) +  m_Q \epsilon \left( \dot{z}_1 q_1 +  \dot{z}_2 q_2 + z_1 \dot{q}_1 + z_2 \dot{q}_2 \right)\\& + m_Q \epsilon^2 ( y_1 \ddot{q}_1 + y_2 \ddot{q}_2),\\
\dot{q}_r =& \omega_r \times q_r, \\
 I_{\eff}\dot{\omega}_r=& q_r \times \left[ u_2 - u_1 +  m_Q L_c (\ddot{q}_2 - \ddot{q}_1)\right.\\ &\qquad\qquad\left. + m_Q \epsilon \left( \dot{z}_2 q_2  - \dot{z}_1 q_1 + z_2 \dot{q}_2 - z_1 \dot{q}_1 \right) + m_Q \epsilon^2 ( y_2 \ddot{q}_2 - y_1 \ddot{q}_1) \right], \\
\epsilon \dot{y}_j =& z_j, \\
\epsilon \dot{z}_j =& \frac1{m_Q}q_j^T \left[  -\bar{c}z_j q_j - \bar{k} y_j q_j - u_j + m_Q (L_c  + \epsilon^2 y_j)\ddot{q}_j + m_Q ( \dot{v}_r + (-1)^j L_r \ddot{q}_r+ ge_3) \right], \\
\dot{q}_j =& \omega_j \times q_j, \\
\dot{\omega}_j =& \frac1{L_c} q_j \times  \left[ \dot{v}_r + (-1)^j L_r \ddot{q}_r + ge_3 - \frac1{m_Q} u_j - \epsilon ( \dot{z}_j q_j + z_j \omega_j) - \epsilon^2 y_j \ddot{q}_j \right], \\
\dot{R}_j =& R_j \hat{\Omega}_j, \\
J_Q \dot{\Omega}_j =& J_Q \Omega_j \times \Omega_j + M_j,  \quad \text{for $j = 1,2$}.
\end{align*} 

The previous system of differential equations can be written as \begin{align}
\dot{x} &= f(t, x, z; \epsilon),\label{eqslow} \\
\epsilon \dot{z} &= g(t, x, z; \epsilon)\label{eqfast},
\end{align} where $f$ and $g$ are smooth functions, $x$ is the vector representing $(x_r, v_r, q_r, \omega_r, q_j, \omega_j, R_j, \Omega_j)$, and $z$ is the vector representing $(y_j, z_j)$, for $j = 1,2$. The above dynamical system is known as a \textit{singular perturbation model} \cite{SPT}, \cite{Khalil}, with \eqref{eqslow} describing the slow dynamics and \eqref{eqfast} describing the fast dynamics.

Evaluating at $\epsilon = 0$, the fast dynamics provide us with algebraic equations that can be solved to obtain $z = h(t, x)$. In particular, \begin{align}
y_j &= \frac1{\bar{k}} q_j^T \left[ -u_j + m_Q L_c\ddot{q}_j + m_Q ( \dot{v}_r + (-1)^j L_r \ddot{q}_r+ ge_3) \right], \label{heq1} \\
z_j &= 0. \label{heq2}
\end{align}

Substituting these equations back into the slow dynamics, we obtain the reduced (slow) model of the control system describing the cooperative task, given by $\dot{x} = f(t, x, h(t, x), 0)$. That is,
\begin{align}
\dot{x}_r &= v_r, \label{xred} \\
m_{\eff} ( \dot{v}_r + ge_3) &= u_1 + u_2 + m_QL_c (\ddot{q}_1 + \ddot{q}_2 ), \label{vred} \\
\dot{q}_r &= \omega_r \times q_r, \label{qrred} \\
 I_{\eff} \omega_r &= q_r \times \left[ u_2 - u_1 + m_Q L_r ( \ddot{q}_{2} - \ddot{q}_1) \right], \label{wrred} \\
 \dot{q}_j &= \omega_j \times q_j, \label{qjred} \\
 m_Q L_c \dot{\omega}_j &= m_Q q_j \times ( \dot{v}_r + (-1)^j L_r \ddot{q}_r+ ge_3) - q_j \times u_j, \label{wjred} \\
\dot{R}_j &= R_j \hat{\Omega}_j, \label{Rkinred} \\
J_Q \dot{\Omega}_j &= J_Q \Omega_j \times \Omega_j + M_j,  \quad \text{for $j = 1,2$} \label{Rdynred}.
\end{align} 
Note that the previous dynamical control system is equivalent to the model with inelastic cables (that is, where $l \equiv L$). Achieving exponentially stable tracking of the reduced model on some set of initial conditions will guarantee exponentially stable tracking in some subset of those initial conditions -- whose relative size depends on $\epsilon$. This claim will be formalized in the subsequent section, and will make direct use of Theorem 11.2 in $\cite{Khalil}$. We will work within this reduced model to design geometric controllers towards the end of tracking the position and attitude of the rod.

\begin{remark}
As discussed in \cite{Elastic}, \cite{SreLeePICDC13}, it is more physically realistic to model our system by a hybrid dynamical system which transitions between the cases of a taut cable (with positive tension magnitude) and a cable with zero tension. However, such a consideration will not play a role in the control design and subsequent analysis, so we will omit such a development.
\end{remark}

\end{section}



\begin{section}{Control Design for Position and Attitude Trajectory Tracking of the Rigid Bar}\label{sec5}

In this section, we will design a thrust controller $u_j\in\mathbb{R}^{3}$ such that the position and attitude of the rigid rod reach a desired position $\tilde{x}_r\in\mathbb{R}^{3}$ and attitude $\tilde{q}_r\in S^{2}$. The strategy is to split $u_j$ into its components which are parallel and perpendicular to the cable attitudes $q_j$. This is motivated by the fact that these components appear independently from one another within the system, and as such can be considered as decoupled controllers. 

Notice that equations \eqref{Rkinred}-\eqref{Rdynred} describing the quadrotor attitude are independent from the rest of the dynamical system, and the moment controllers $M_j$ appear exclusively within them. Moreover, these equations are identical to those that appear in \cite{LLM}, \cite{Lee-TAC}, for which $M_j$ was designed to attain almost-global exponential stability. We use the same controller for the attitude of the quadrotors, and disregard the equations for the remainder of this paper.

We will also introduce configuration error functions for each state variable, from which we may derive our error dynamics. The controls will be selected such that the origin is an exponentially stable equilibrium point of the error dynamics. In particular, the controls will take the form of a feedback linearization together with a PD controller. 

\subsection{Error dynamics}

We begin by further simplifying the dynamical system \eqref{xred}-\eqref{wjred}. In particular, we find an equation for $\ddot{q}_j$ that we will substitute into equations \eqref{vred} and \eqref{wjred}. By differentiating \eqref{qjred} and expanding it with the vector triple product identity, it can be shown that $\ddot{q}_j = \dot{\omega}_j \times q_j - ||\omega_j ||^2 q_j$. Now we may substitute \eqref{wjred} in for $\dot{\omega}_j$ to find that \begin{align*}
m_Q L_c \ddot{q}_j &= -q_j \times (m_Q L_c\dot{\omega}_j) - m_Q L_c||\omega_j||^2 q_j \\
&= -q_j \times \left[ m_Q q_j \times ( \dot{v}_r + (-1)^j L_r \ddot{q}_r + ge_3) - q_j \times u_j \right] - m_Q L_c ||\omega_j||^2 q_j\\
&= -m_Q \hat{q}_j^2  \left[ \dot{v}_r + (-1)^j L_r \ddot{q}_r + ge_3\right] + \hat{q}_j^2 u_j - m_Q L_c ||\omega_j||^2 q_j \\ 
&= m_Q (\dot{v}_r^{\perp_j} + (-1)^j L_r \ddot{q}_r^{\perp_j} + ge_3^{\perp_j}) - u_j^{\perp_j} - m_Q L_c ||\omega_j||^2 q_j,
\end{align*}where ``$\perp_j$'' stands for the component of the vector that is perpendicular to the cable attitude $q_j\in S^{2}$. Similarly, for every vector $v \in \R^3$, we define the component $v^{\parallel_j}$ parallel to $q_j$ such that $v = v^{\parallel_j} + v^{\perp_j}$. Notice that $v^{\parallel_j} = (v^T q_j)q_j$ and $v^{\perp_j} = -\hat{q}_j^2 v$. Substituting the above equation for $m_Q L_c \ddot{q}_j$ into \eqref{vred} and making use of the fact that $m_{\eff} = 2m_Q + m_r$, we obtain \begin{align*}
m_{r} ( \dot{v}_r + ge_3) &= \left[ u_1 - m_Q( \dot{v}_r + ge_3) + m_Q L_c \ddot{q}_1 \right] + \left[ u_2 - m_Q(\dot{v}_r + ge_3)+ m_QL_c\ddot{q}_2 \right],\\
&= \left[ u_1^{\parallel_1} - m_Q (\dot{v}_r^{\parallel_1} + ge_3^{\parallel_1} + L_c ||\omega_1||^2 q_1 - L_r \ddot{q}_r^{\perp_1}) \right]\\& \quad + \left[ u_2^{\parallel_2} - m_Q (\dot{v}_r^{\parallel_2} + ge_3^{\parallel_2} + L_c ||\omega_2||^2 q_2 + L_r \ddot{q}_r^{\perp_2}) \right].
\end{align*}

Define the parallel component of the thrust controllers as \begin{equation}\label{uparallel}u_j^{\parallel_j} := \mu_j + m_Q (\ddot{x}_r^{\parallel_j} + ge_3^{\parallel_j} + L_c ||\omega_j||^2 q_j + (-1)^j L_r \ddot{q}_r^{\parallel_j}),\end{equation} where $\mu_j$ is an additional control to be designed later -- and which is constrained to be parallel to $q_j$. Repeating this procedure with \eqref{wrred} yields 
\begin{align*}
 I_{\eff} \omega_r &= q_r \times \left[ u_2 + m_Q L_c \ddot{q}_2 \right]  - q_r \times \left[ u_1 + m_Q L_c \ddot{q}_1 \right], \\
&= q_r \times \left[ u_2^{\parallel_2} + m_Q (\dot{v}_r^{\perp_2} + ge_3^{\perp_2} - L_c ||\omega_2||^2 q_2 +  L_r \ddot{q}_r^{\perp_2}) \right] \\ &- q_r \times \left[ u_1^{\parallel_1} + m_Q (\dot{v}_r^{\perp_1} + ge_3^{\perp_1} - L_c ||\omega_1||^2 q_1  - L_r \ddot{q}_r^{\perp_1}) \right]\\
&= q_r \times \left[ \mu_2 - \mu_1 + 2m_Q L_r \ddot{q}_r \right].
\end{align*}

Further making use of the fact that $I_{\eff} = 2m_QL_r + \frac23 m_r L_r$, we obtain the following dynamical system
\begin{align}
\dot{x}_r &= v_r, \label{xkin2}\\
m_r ( \dot{v}_r + ge_3) &= \mu_1 + \mu_2, \label{xdyn2}\\
\dot{q}_r &= \omega_r \times q_r, \label{qrkin2}\\
\frac23 m_r L_r \dot{\omega}_r &= q_r \times ( \mu_2 - \mu_1 ), \label{qrdyn2} \\
\dot{q}_j &= \omega_j \times q_j, \label{qjkin2} \\
m_Q L_c \dot{\omega}_j &= m_Q q_j \times ( \ddot{x}_r + (-1)^j L_r \ddot{q}_r + ge_3) - q_j \times u_j. \label{qjdyn2}
\end{align}

Equations \eqref{xkin2}, \eqref{qrkin2}, and \eqref{qjkin2} are the kinematic equations for the load position, load attitude, and cable attitude, respectively, while equations \eqref{xdyn2}, \eqref{qrdyn2}, and \eqref{qjdyn2} describe the dynamics. Notice also that the cross product $q_j \times u_j$ in equation \eqref{qjdyn2} will annihilate $u_j^{\parallel_j}$, so that we may view $\mu_1, \mu_2,$ and $q_j \times u_j$ as completely independent. We now introduce the tracking errors for the position and velocity of the rod as $e_{x_r} := x_{r} - \tilde{x}_{r}$ and $e_{v_r} := v_r - \tilde{v}_r = \dot{e}_x$, respectively. These clearly have the desirable property that $e_{x_r} = e_{v_r} = 0$ implies $x_r = \tilde{x}_r \text{ and } v_r = \tilde{v}_r$. 

\vspace{.2cm}

We define the tracking error for the attitude of the rod as $e_{q_r}:= \tilde{q}_r \times q_r$, which is a tangent vector in $T_{q_r} S^2$ with the property that if $e_{q_r} = 0$ then $q_r = \tilde{q}_r \text{ or } q_r = -\tilde{q}_r$. Note that $\omega_{r}$ and $\tilde{\omega}_r$ will in general belong to different tangent spaces, and so to compare them, we must first translate $\tilde{\omega}_r$ to $T_{q_r}S^2$. After which, they are vectors belonging to the same vector space and can be compared with their difference just as we did with $e_{x_r}$ and $e_{v_r}$. That is, we can define the tracking error for the angular velocity of the rod as $e_{\omega_r} = \omega_r + \hat{q}_r^2 \tilde{\omega}_r.$ This analysis is not specific to $q_r$ or $\omega_r$, but instead relies only on the manifold structure of $S^2$. Hence we may use the same configuration errors to compare $q_j, \omega_j$ with their respective desired trajectories for $j = 1,2$.

Similarly, we define the configuration error function on $S^2$ for $q_j$ as $\Psi_{j} = 1 - \tilde{q}_{j}^T q_j$, for $j = 1,2,r$, which clearly satisfies $0 \le \Psi_{j} \le 2$, with $\Psi_j = 0$ if and only if $q_j = \tilde{q}_j$. Observe that for some $\theta \in [0, 2 \pi)$, we may write $\Psi_{j} = 1 - \cos(\theta)$ and $||e_{q_r}||^2 = \sin^2(\theta)$. With this representation, is it easy to see that whenever $\Psi_{j} \le \psi_{j} < 1$, we have $\frac12 ||e_{q_j}||^2 \le \Psi_{j} \le \frac1{2 - \psi_{j}} ||e_{q_j}||^2$. 
 
Note that, from \eqref{xdyn2}, the error dynamics for the load position are given by $m_r(\dot{e}_{v_r} + \dot{\tilde{v}}_r + ge_3) = \mu_1 + \mu_2$.

\noindent while the error dynamics for the load attitude are described by \begin{align*}
\frac23 m_r L_r\dot{e}_{\omega_r} &= \frac23 m_r L_r \left[ \dot{\omega}_r + \dot{q}_r \times (q_r \times \tilde{\omega}_r) + q_r \times (\dot{q}_r \times \tilde{\omega}_r) + \hat{q}_r^2 \dot{\tilde{\omega}}_r \right] \\
&= \frac23 m_r L_r \left[\dot{\omega}_r + \hat{q}_r^2 \dot{\tilde{\omega}}_r + (\dot{q}_r^T \tilde{\omega}_r)q_r + (q_r^T \tilde{\omega}_r)\dot{q}_r\right] \\
&= \hat{q}_r \left[ \mu_2 - \mu_1 + \frac23 m_r L_r \hat{q}_r \dot{\tilde{\omega}}_r -  \frac23 m_r L_r(q_r^T \tilde{\omega}_r) (e_{\omega_r} - \hat{q}_r^2 \tilde{\omega}_r) \right] + \frac23 m_r L_r (\dot{q}_r^T \tilde{\omega}_r)q_r.
\end{align*}

This motivates us to choose our "desired" controls $\tilde{\mu}_1, \tilde{\mu}_2$ such that their sum is a feedback linearization plus a PD controller for the above error dynamics. That is \begin{align}
\tilde{\mu}_1 + \tilde{\mu}_2 &= m_r(\dot{\tilde{v}}_r + ge_3 - k_{v_r} e_{v_r} - k_{x_r} e_{x_r}), \label{mutilde}\\
\tilde{\mu}_2 - \tilde{\mu}_1 &= \frac23 m_r L_r \left[ -\hat{q}_r \dot{\tilde{\omega}}_r +  (q_r^T \tilde{\omega}_r) \hat{q}_r^2 \tilde{\omega}_r  -  (q_r^T \tilde{\omega}_r)e_{\omega_r} +  k_{\omega_r} \hat{q}_r e_{\omega_r} + k_{q_r}\hat{q}_r e_{q_r} \right] \label{mutildeq}.
\end{align}


Note, however, that we can not simply take $\mu_j = \tilde{\mu}_j$, as $\tilde{\mu}_j$ is not guaranteed to satisfy the constraint of parallelism to $q_j$. Therefore, we chose $\mu_j = (I + \hat{q}_j^2)\tilde{\mu}_j$ and elect our desired cable attitudes $\tilde{q}_j$ such that $\hat{\tilde{q}}_j^2 \tilde{\mu}_j = 0$. In particular, we define $\displaystyle{\tilde{q}_j = -\frac{\tilde{\mu}_j}{||\mu_j||}}$. With such a choice, the error dynamics of the translational and angular velocities of the rod become \begin{align*}
\dot{e}_{v_r} &= -k_{v_r} e_{v_r} - k_{x_r} e_{x_r} + \frac1{m_r} ||\tilde{\mu}_1||\hat{q}_1 e_{q_1} + \frac1{m_r} ||\tilde{\mu}_2|| \hat{q}_2 e_{q_2}, \\
\dot{e}_{\omega_r} &= (\dot{q}_r^T \tilde{\omega}_r)q_r - k_{\omega_r} e_{\omega_r} - k_{q_r} e_{q_r} + \frac3{2 m_r L_r}\hat{q}_r (||\tilde{\mu}_2|| \hat{q}_2 e_{q_2} - ||\tilde{\mu}_1|| \hat{q}_1 e_{q_1}).
\end{align*}

Next, we design the perpendicular component $u_j^{\perp}$ of the thrust controller in such a way $q_j$ approaches $\tilde{q}_j$ asymptotically. Recall that the dynamics for the cables are \begin{align*}
\dot{q}_j &= \omega_j \times q_j, \\
m_Q L_c \dot{\omega}_j &= m_Q q_j \times ( \ddot{x}_r + (-1)^j L_r \ddot{q}_r + ge_3) - q_j \times u_j.
\end{align*}

Repeating the procedure that we used to find $e_{\omega_r}$, we find that the error dynamics for the angular velocities of the cables as \begin{align*}
m_Q L_c \dot{e}_{\omega_j} &= m_Q L_c \dot{\omega}_j + m_Q L_c \left[ \dot{q}_j \times (q_j \times \tilde{\omega}_j) + q_j \times (\dot{q}_j \times \tilde{\omega}_j) + \hat{q}_j^2 \dot{\tilde{\omega}}_j) \right] \\
&= m_Q \hat{q}_j ( \ddot{x}_r + (-1)^j L_r \ddot{q}_r + ge_3) + m_Q L_c \left[\hat{q}_j^2 \dot{\tilde{\omega}}_j - (q_j^T \tilde{\omega}_j)\hat{q}_j \omega_j + (\dot{q}_j^T \tilde{\omega}_j)q_j\right] - \hat{q}_j u_j \\
&=m_Q \hat{q}_j \left[ \ddot{x}_r + (-1)^j L_r \ddot{q}_r + ge_3 + L_c \hat{q}_j \dot{\tilde{\omega}}_j - L_c(q_j^T \tilde{\omega}_j) \omega_j \right] + m_Q L_c(\dot{q}_j^T \tilde{\omega}_j)q_j - \hat{q}_j u_j. 
\end{align*}

Since $\hat{q}_j u_j = \hat{q}_j u_j^{\perp_j}$, we can choose $u_j^{\perp_j}$ such that \begin{align}
\hat{q}_j u_j^{\perp_j} &= m_Q \hat{q}_j \left[ \ddot{x}_r + (-1)^j L_r \ddot{q}_r + ge_3 + L_c \hat{q}_j \dot{\tilde{\omega}}_j - L_c(q_j^T \tilde{\omega}_j) \omega_j \right] + m_Q L_c \left[ k_{q_j} e_{q_j} + k_{\omega_j} e_{\omega_j}\right] \nonumber\\
u_j^{\perp_j} &= -m_Q \hat{q}_j^2\left[ \ddot{x}_r + (-1)^j L_r \ddot{q}_r + ge_3 + L_c \hat{q}_j \dot{\tilde{\omega}}_j - L_c(q_j^T \tilde{\omega}_j) \omega_j \right] - m_Q L_c \hat{q}_j \left[ k_{q_j} e_{q_j}  + k_{\omega_j} e_{\omega_j} \right]\label{ujperp},
\end{align} 
which clearly satisfies the orthogonality constraints and leads to the following error dynamics for $\omega_j$ $$\dot{e}_{\omega_j} =  (\dot{q}_j^T \tilde{\omega}_j)q_j - k_{q_j} e_{q_j} - k_{\omega_j} e_{\omega_j}.$$ 

\vspace{.2cm}
\subsection{Control design for the reduced model}
We now show that an appropriate choice of gains will cause the origin of the error dynamics for the reduced model to be exponentially stable, where the control inputs are defined as above. In particular, we use a Lyapunov candidate to find sufficient conditions under which the origin of the error dynamics is exponentially stable. These conditions are written in terms of the gains and some arbitrary Lyapunov parameters. Then, we show that these conditions can be satisfied, provided that the initial errors in the cable attitudes are sufficiently small.

Consider the indexing set $\mathcal{I} = {1, 2, r}$, the domain \begin{align*}\mathcal{D} = \{ (e_{x_r}, e_{v_r}, e_{q_r}, e_{\omega_r}, e_{q_1},& e_{\omega_1}, e_{q_2}, e_{\omega_2}) | \  ||e_{x_r}|| \le \bar{e}_{x_r}, ||e_{v_r}|| \le \bar{e}_{v_r}, e_{\omega_j} \le \bar{e}_{\omega_j}, \ \Psi_{j} \le \psi_j\}, \end{align*} and define a Lyapunov candidate function on $\mathcal{D}$ given by
\begin{align*}
V = \frac12 ||e_{v_r}||^2 + \frac12 k_{x_r} ||e_{x_r}||^2 + c_{x_r} e_{x_r}^T e_{v_r} + \sum_{j\in\mathcal{I}} \left[ \frac12 ||e_{\omega_j}||^2 + k_{q_j} \Psi_{j} + c_{q_j} e_{q_j}^T e_{\omega_j} \right],
\end{align*}where $c_{x_r}$ and $c_{q_j}$ are positive real numbers for $j \in \mathcal{I}$. Observe that $V_x$ defined as $V_{x} := \frac12 ||e_{v_{r}}||^2 + c_{x_r} e_{x_r}^T e_{v_r}+ \frac12 k_{x_{r}} ||e_{x_r}||^2$ can be bounded from above and below as $\frac12 z_{x}^T \underbar{P}_{x} z_{x} \le \ V_{x} \le \frac12 z_{x}^T \bar{P}_{x} z_{x}$ where $\underbar{P}_x = \begin{bmatrix}
k_{x_r} & -c_{x_r} \\  
-c_{x_r} & 1
\end{bmatrix}, \bar{P}_x = \begin{bmatrix}
k_{x_r} & c_{x_r} \\  
c_{x_r} & 1
\end{bmatrix}$, and $z_x = \begin{bmatrix} ||e_{x_r}|| & ||e_{v_r}|| \end{bmatrix}^T.$ Further note that both $\underbar{P}_x$ and $\bar{P}_x$ are positive-definite provided that $c_{x_r} < \sqrt{k_{x_r}}$.

Similarly, for $j\in\mathcal{I}$, we define $V_{q_j} := \frac12 ||e_{\omega_j}||^2 + c_{q_j} e_{q_j}^T e_{\omega_j} + k_{q_j} \Psi_{q_j}$, which is bounded as $\frac12 z_{q_j}^T \underbar{P}_{q_j} z_{q_j} \le V_{q_j} \le \frac12 z_{q_j}^T \bar{P}_{q_j} z_{q_j},$ where $\underbar{P}_{q_j} = \begin{bmatrix}
k_{q_j} & -c_{q_j} \\  
-c_{q_j} & 1\end{bmatrix}, \bar{P}_{q_j} = \begin{bmatrix}
\frac{2k_{q_j}}{2 - \psi_{q_j}} & c_{q_j} \\  
c_{q_j} & 1\end{bmatrix},$ and $z_{q_j} = \begin{bmatrix} ||e_{q_j}|| & ||e_{\omega_j}|| \end{bmatrix}^T$. As before,  $\underbar{P}_{q_j}$ and $\bar{P}_{q_j}$ are positive-definite  when $c_{q_j} < \sqrt{k_{q_j}}.$

Observing that $V = V_x + \sum_{j \in \mathcal{I}} V_{q_j}$, we then have that our Lyapunov candidate is bounded as $\frac12 z^T \underbar{P} z \le V \le \frac12 z^T \bar{P} z,$ where $\footnotesize{z = \begin{bmatrix} ||e_{x_r}|| & ||e_{v_r}|| & ||e_{q_r}|| & ||e_{\omega_r}|| & ||e_{q_1}|| & ||e_{\omega_1}|| & ||e_{q_2}|| & ||e_{\omega_2}|| \end{bmatrix}^T},$

 \begin{align*}
\underbar{P} &= \begin{bmatrix} \underbar{P}_x & 0 & 0 & 0 \\ 0 & \underbar{P}_{q_r} & 0 & 0 \\ 0 & 0 & \underbar{P}_{q_1} & 0 \\ 0 & 0 & 0 & \underbar{P}_{q_2}  \end{bmatrix}, \quad \text{ and } \quad  \bar{P} = \begin{bmatrix} \bar{P}_x & 0 & 0 & 0 \\ 0 & \bar{P}_{q_r} & 0 & 0 \\ 0 & 0 & \bar{P}_{q_1} & 0 \\ 0 & 0 & 0 & \bar{P}_{q_2}  \end{bmatrix}. 
\end{align*} where $\underbar{P}, \bar{P}$ are positive-definite for $c_{x_r} < \sqrt{k_{x_r}}$ and $c_{q_j} < \sqrt{k_{q_j}}$. Next, note that by the invariance of circular shifts of the scalar triple product and the fact that $q_j^T e_{q_j} = 0$, we have:
\begin{align*}
\frac{d}{dt}\Psi_{q_j} &= -\tilde{q}^T_j \dot{q}_j - q_j^T \dot{\tilde{q}}_j = -\tilde{q}_j^T (\omega_j \times q_j) - q_j^T (\tilde{\omega}_j \times \tilde{q}_j)\\ 
&= \omega_j^T (\tilde{q}_j \times q_j) - \tilde{\omega}_j^T (\tilde{q}_j \times q_j) \\
&= (\omega_j - \tilde{\omega}_j)^T e_{q_j} = (\omega_j + (q_j^T \tilde{\omega}) q_j - \tilde{\omega}_j)^T e_{q_j} = e_{\omega_j}^T e_{q_j}.
\end{align*}
Additionally, from the vector triple product, we see
\begin{align*}
\dot{e}_{q_j} &= (\dot{\tilde{q}}_j \times q_j) + (\tilde{q}_j \times \dot{q}_j) = (\tilde{\omega}_j \times \tilde{q}_j) \times q_j - (\omega_j \times q_j) \times \tilde{q}_j \\
&= \tilde{\omega}_j \times (\tilde{q}_j \times q_j) - \tilde{q}_j \times (\tilde{\omega}_j \times q_j) - \omega_j \times (q_j \times \tilde{q}_j) + q_j \times (\omega_j \times \tilde{q}_j) \\
&= (\omega_j + \tilde{\omega}_j) \times e_{q_j} + (\tilde{q}_j^T q_j)e_{\omega_j} -  (\tilde{q}_j^T q_j)  (q_j^T \tilde{\omega}_j)q_j \\
&= e_{\omega_j} \times e_{q_j} + (\tilde{q}_j^T q_j)e_{\omega_j} + 2\tilde{\omega}_j \times e_{q_j} -  (\tilde{q}_j^T q_j)  (q_j^T \tilde{\omega}_j)q_j,
\end{align*}
so that
\begin{align*}
\dot{e}_{q_j}^T e_{\omega_j} &= (\tilde{q}_j^T q_j)||e_{\omega_j}||^2 + 2(\tilde{\omega}_j \times e_{q_j})^T e_{\omega_j} \le ||e_{\omega_j}||^2 + C_{q_j} ||e_{q_j}|| ||e_{\omega_j}||,
\end{align*} where $C_{q_j} \le 2\sup ||\tilde{\omega}_j||$ is a non-negative constant. 
Therefore, the time derivative of the proposed Lyapunov function is bounded as
 \begin{align*}
\dot{V}  \le& -(k_{v_r} - c_{x_r}) ||e_{v_r}||^2 + c_{x_r} k_{v_r} ||e_{v_r}|| ||e_{x_r}|| - c_{x_r}k_{x_r}||e_{x_r}||^2 \\&+ \frac1{m_r}(||e_{v_r}|| + c_{x_r} ||e_{x_r}||) ||Y|| 
  -(k_{\omega_r} - c_{q_r})||e_{\omega_r}||^2\\& + c_{q_r}(k_{\omega_r} + C_{q_r})||e_{\omega_r}||||e_{q_r}|| - c_{q_r} k_{q_r} ||e_{q_r}||^2 - c_{q_j} k_{q_j} ||e_{q_j}||^2\\ &+ \frac2{3 m_r L_r} (||e_{\omega_r}|| + c_{q_r} ||e_{q_r}||) ||Y||   -(k_{\omega_j} - c_{q_j}) ||e_{\omega_j}||^2 \\&+ c_{q_j} (k_{\omega_j} + C_{q_j} ) ||e_{q_j}|| ||e_{\omega_j}||,
\end{align*} where $Y$ satisfies the inequality \begin{align*}
||Y|| \le \left[ m_r (k_{v_r} ||e_{v_r}|| + k_{x_r} ||e_{x_r}||) + \frac23 m_r L_r ((C^2_{q_r} + k_{\omega_r}) ||e_{\omega_r}|| + k_{q_r}||e_{q_r}||) + C\right] (||e_{q_1}|| + ||e_{q_2}||)
\end{align*} for a non-negative constant $C \le m_r \sup ||\dot{\tilde{v}}_r|| + \frac23 m_r L_r \sup ||\dot{\tilde{\omega}}_r||$. Furthermore, within $\mathcal{D}$, we have $||e_{q_j}|| \le \sqrt{\psi_{q_j} (2 - \psi_{q_j})} := \alpha_j$. With $\alpha = 2\max \{ \alpha_1, \alpha_2 \}$ and $I_r = \frac23 m_r L_r$, we have \begin{align*}
 \frac1{m_r}(||e_{v_r}|| + c_{x_r} ||e_{x_r}||) ||Y|| \le & \alpha(k_{v_r}||e_{v_r}||^2 + c_{x_r}k_{v_r}||e_{x_r}||||e_{v_r}|| + c_{x_r} ||e_{x_r}||^2 ) \\ 
 & + \frac{I_r}{m_r} \alpha (||e_{v_r}|| + c_{x_r} ||e_{x_r}||)((C_{q_r}^2 + k_{\omega_r})||e_{\omega_r}|| + k_{q_r}||e_{q_r}||) \\
 & + \frac{1}{m_r}((C + m_r k_{x_r} \bar{e}_{x_r})||e_{v_r}|| + c_{x_r} C ||e_{x_r}||) (||e_{q_1}|| + ||e_{q_2}||),\end{align*}  and, 
\begin{align*} (||e_{\omega_r}|| + c_{q_r} ||e_{q_r}||) ||Y|| \le & I_r\alpha ((C_{q_r}^2 + k_{\omega_r})||e_{\omega_r}||^2\\& + c_{q_r}(C_{q_r}^2 + k_{\omega_r})||e_{q_r}||||e_{\omega_r}|| + c_{q_r} k_{q_r}||e_{q_r}||^2) \\
 & + m_r \alpha ( ||e_{\omega_r}|| + c_{q_r} ||e_{q_r}|| )(k_{v_r} ||e_{v_r}|| + k_{x_r} ||e_{x_r}||) \\
 & + ((C + I_r k_{q_r} \alpha_r) ||e_{\omega_r}|| + c_{q_r}C ||e_{q_r}|| )(||e_{q_1}|| + ||e_{q_2}||).
 \end{align*}Applying this inequality directly to the above bound for $\dot{V}$, we find that $\dot{V} \le -z^T \mathcal{W} z$, where $\mathcal{W} = \mathcal{W}_1 + \mathcal{W}_2$, and for $j = 1,2$,  $\mathcal{W}_j$ is the $6\times6$ matrix defined as $$\mathcal{W}_j = \begin{bmatrix}
 W_{x_r} & -\frac12 W_{x_r, q_r} & -\frac12 W_{x_r, q_j} \\
 -\frac12 W_{x_r, q_r} & W_{q_r} & -\frac12 W_{q_r, q_j} \\
-\frac12 W_{x_r, q_j} & -\frac12 W_{q_r, q_j} & W_{q_j} \\
 \end{bmatrix}$$
 for $2\times2$ sub-matrices given by
 \begin{align*}
 &W_{x_r} = 
 \begin{bmatrix}
 c_{x_r} (1 - \alpha) k_{x_r} & -\frac12 (1 + \alpha) c_{x_r} k_{v_r} \\
 -\frac12 (1 + \alpha) c_{x_r} k_{v_r} & (1 - \alpha) k_{v_r} - c_{x_r}
 \end{bmatrix}, \\
 &W_{q_j} = 
 \begin{bmatrix}
c_{q_j} k_{q_j}  & -\frac12 c_{q_j}(k_{\omega_j} + C_{q_j}) \\
 -\frac12 c_{q_j}(k_{\omega_j} + C_{q_j})  & k_{\omega_j} - c_{q_j}
 \end{bmatrix}, \\
 &W_{q_r} = \begin{bmatrix}
 (1 - \alpha)c_{q_r} k_{q_r} & - \frac12 c_{q_r}((1 + \alpha)k_{\omega_r} + C_{q_r} + \alpha C_{q_r}^2) \\
 -  \frac12 c_{q_r}((1 + \alpha)k_{\omega_r} + C_{q_r} + \alpha C_{q_r}^2) &  (1 - \alpha)k_{\omega_r} - c_{q_r} - \alpha C_{q_r}^2
 \end{bmatrix}\\
& W_{x_r, q_r} = \alpha
\begin{bmatrix}
\frac{I_r}{m_r}c_{x_r}k_{q_r} + \frac{m_r}{I_r} c_{q_r} k_{x_r} & \frac{I_r}{m_r}c_{x_r}(C_{q_r}^2 + k_{\omega_r}) + \frac{m_r}{I_r} k_{x_r}  \\
\frac{I_r}{m_r} k_{q_r} + \frac{m_r}{I_r} c_{q_r} k_{v_r} & \frac{I_r}{m_r}(C_{q_r}^2 + k_{\omega_r}) + \frac{m_r}{I_r} k_{v_r}
\end{bmatrix}, \\
&W_{x_r, q_j} = \frac1{m_r}
\begin{bmatrix}
c_{x_r} C & 0 \\
C + m_r k_{x_r} \bar{e}_{x_r} & 0
\end{bmatrix}, \,\,
W_{q_r, q_j} = \frac1{I_r}
\begin{bmatrix}
c_{q_r} C & 0 \\
C + I_r k_{q_r} \alpha_r & 0
\end{bmatrix}.
\end{align*} Note that $\mathcal{W}$ is not necessarily symmetric, so the quadratic form defined by $z^T \mathcal{W} z$ is positive definite if and only if the symmetric part of $\mathcal{W}$—that is $\frac12 (\mathcal{W} + \mathcal{W}^T)$—is a positive definite matrix. The following Theorem states, in essence, that the gains and Lyapunov constants can be chosen such that $\underbar{P}, \bar{P},$ and $\frac12 (\mathcal{W} + \mathcal{W}^T)$ are simultaneously positive-definite—thus ensuring that the origin is exponentially stable—provided that the initial errors in the cable attitudes are sufficiently small. We now state the main result of this work:

\begin{theorem}\label{gainsth}
Consider the control system with disturbances defined by equations \ref{xkin2} - \ref{qjdyn2} with control inputs \ref{mutilde}, \ref{mutildeq}, and \ref{ujperp}. For sufficiently small $\alpha$, there exists control gains $k_{x_r}$, $k_{v_r}$, $k_{q_r}$, $k_{\omega_r}$, $k_{q_1}$, $k_{\omega_1}$, $k_{q_2}$, and $k_{\omega_2}$ such that the zero equilibrium of the tracking errors $e_{x_r}$, $e_{v_r}$, $e_{q_r}$, $e_{\omega_r}$, $e_{q_1}$, $e_{\omega_1}$, $e_{q_2}$ and $e_{\omega_2}$ is exponentially stable.
\end{theorem}

Before proving this theorem, we state some general facts about positive-definite matrices that will be used multiple times in the proof.

\begin{lemma}{(Facts about positive-definite matrices)}\label{posdeflemma}
Let $A \succ 0$ and $B \succeq 0$ with $\lambda_{\min}(A) > \lambda_{\max}(B)$, and $M$ arbitrary, be $n\times n$ matrices. Then, the following characterizations hold
\begin{enumerate}
\item $||x||||y|| \lambda_{\min}(A) \le x^T A y \le ||x||||y|| \lambda_{\max}(A)$,
\item $M^T A M \succeq 0$,
\item $A - B \succ 0$ with $\lambda_{\min}(A-B) \ge \lambda_{\min}(A) - \lambda_{\max}(B)$ and $\lambda_{\max}(A-B) \le \lambda_{\max}(A)$.
\end{enumerate}
\end{lemma}

\subsection*{Proof of Theorem \ref{gainsth}} 
Denote the symmetric part of $\mathcal{W}$ by $\bar{\mathcal{W}} = \frac12 (\mathcal{W} + \mathcal{W}^T)$ and similarly define the symmetric parts of the submatrices by $\bar{W}_{x_r, q_r}, \bar{W}_{x_r, q_j},$ and $\bar{W}_{q_r. q_j}$. It is clear that $\bar{\mathcal{W}}$ can be expressed in the form $\bar{\mathcal{W}} = \begin{bmatrix} P & S \\ S^T & Q \end{bmatrix}$, where $P = \begin{bmatrix} W_{x_r} & -\frac12 \bar{W}_{x_r, q_r} \\ -\frac12 \bar{W}_{x_r, q_r} & {W}_{q_r} \end{bmatrix}$, $S = -\frac12 \begin{bmatrix} \bar{W}_{x_r, q_j} & \bar{W}_{q_r, q_j} \end{bmatrix}^T$, and $Q = {W}_{q_j}$. 

Now, observe that $\mathcal{W}$ can be decomposed as:
\begin{align}\label{Schur}
\begin{bmatrix} P & S \\ S^T & Q \end{bmatrix} = \begin{bmatrix} I & SQ^{-1} \\ 0 & I \end{bmatrix} \begin{bmatrix} P - SQ^{-1}S^T  & 0 \\ 0 & Q \end{bmatrix} \begin{bmatrix} I & SQ^{-1} \\ 0 & I \end{bmatrix}^T
\end{align}
Where $P - SQ^{-1}S^T$ is often referred to as the \textit{Schur complement} of $Q$. From Lemma \ref{posdeflemma}, it then follows that $\mathcal{\bar{W}} \succ 0$ if and only if $P - SQ^{-1} S^T \succ 0$ and $Q \succ 0.$ Note that $P - SQ^{-1}S^T$ can itself be expressed in form of a $4\times4$ block matrix given by


$$P - SQ^{-1}S^T = \begin{bmatrix}
{W}_{x_r} - \frac14 \bar{W}_{x_r, q_j} {W}_{q_j}^{-1} \bar{W}_{x_r, q_j} & -\frac12 \bar{W}_{x_r,q_r} - \frac14 \bar{W}_{q_r,q_j}{W}_{q_j}^{-1}\bar{W}_{x_r,q_j} \\ -\frac12 \bar{W}_{x_r,q_r} - \frac14 \bar{W}_{q_r,q_j}{W}_{q_j}^{-1}\bar{W}_{x_r,q_j} & {W}_{q_r} -\frac14 \bar{W}_{q_r,q_j}{W}_{q_j}^{-1}\bar{W}_{q_r,q_j}
\end{bmatrix}.$$

Repeating the previous analysis, but now on $P - SQ^{-1}S^T$, we find that $\bar{\mathcal{W}} \succ 0$ if and only if the following three conditions hold: \begin{enumerate}
\item[(1)] $W_{q_j} \succ 0$, \quad (2) $W_{q_r} - \frac14 \bar{W}_{q_r, q_j} W_{q_j}^{-1} \bar{W}_{q_r, q_j} \succ 0$,
\item [(3)]  \begin{align*}
0\prec W_{x_r} &- \frac14 \bar{W}_{x_r, q_j} W_{q_j}^{-1} \bar{W}_{x_r, q_j} - \frac14 \bar{W}_{x_r, q_r} (W_{q_r} - \frac14 \bar{W}_{q_r, q_j} W_{q_j}^{-1} \bar{W}_{q_r,q_j})^{-1} \bar{W}_{x_r, q_r} \\
& - \frac18 \bar{W}_{x_r, q_r} (W_{q_r} - \frac14 \bar{W}_{q_r, q_j} W_{q_j}^{-1} \bar{W}_{q_r,q_j})^{-1} \bar{W}_{q_r,q_j}W_{q_j}^{-1}\bar{W}_{x_r, q_j} \\
& - \frac18 \bar{W}_{x_r,q_j} W_{q_j}^{-1} \bar{W}_{q_r,q_j} (W_{q_r} - \frac14 \bar{W}_{q_r, q_j} W_{q_j}^{-1} \bar{W}_{q_r,q_j})^{-1} \bar{W}_{x_r, q_r} \\
& - \frac1{16} \bar{W}_{x_r,q_j} W_{q_j}^{-1} \bar{W}_{q_r,q_j} (W_{q_r} - \frac14 \bar{W}_{q_r, q_j} W_{q_j}^{-1} \bar{W}_{q_r,q_j})^{-1} \bar{W}_{q_r,q_j}W_{q_j}^{-1}\bar{W}_{x_r, q_j}.
\end{align*}
\end{enumerate}

Moreover, the minimum (maximum) eigenvalue of $\mathcal{W}$ is exactly the smallest (largest) of the minimum (maximum) eigenvalues of the three matrices in the above conditions. We now seek to verify that appropriate choices in the gains and constants can be made to satisfy the above conditions. First, by looking at the characteristic equation of $W_{q_j}$, it can be seen that
 $$2\lambda_{\min}(W_{q_j}) = (k_{\omega_j} - c_{q_j} + c_{q_j}k_{q_j}) - \sqrt{(k_{\omega_j} - c_{q_j} - c_{q_j}k_{q_j}) + c_{q_j}^2(k_{\omega_j} + C_{q_j})^2 }.$$

By taking $c_{q_j}$ sufficiently small, $k_{\omega_j}$ sufficiently large, and $k_{q_j} = \frac{k_{\omega_j}}{c_{q_j}}$, $\lambda_{\min}(W_{q_j})$ can be made arbitrarily large. Consequently, $\lambda_{\max}(W_{q_j}^{-1}) = \lambda_{\min}(W_{q_j})^{-1}$ can be made arbitrarily small (and positive). From Lemma \ref{posdeflemma}, we have that $\frac14 \bar{W}_{q_r, q_j} W_{q_j}^{-1} \bar{W}_{q_r, q_j} \succeq 0$ and, since $\bar{W}_{q_r, q_j}$ and $W_{q_j}$ are independent, we can shrink the maximum eigenvalue of $\bar{W}_{q_r, q_j} W_{q_j}^{-1} \bar{W}_{q_r, q_j}$ arbitrarily by shrinking the maximum eigenvalue of $W_{q_j}^{-1}$.

Similarly, from the characteristic equation of $W_{q_r}$, we see that the eigenvalues satisfy \begin{align*}2 \lambda_{\pm} &= ((1 - \alpha)c_{q_r}k_{q_r} + (1 - \alpha)k_{\omega_r} - c_{q_r} - \alpha C_{q_r}^2)\\&\pm  \sqrt{((1 - \alpha)c_{q_r}k_{q_r} - (1 - \alpha)k_{\omega_r} + c_{q_r} + \alpha C_{q_r}^2)^2 + c_{q_r}^2((1 + \alpha)k_{\omega_r} + C_{q_r} + \alpha C_{q_r}^2)^2}.\end{align*} Now choose $k_{\omega_r} = \frac{\alpha}{1 - \alpha}C_{q_r}^2 + c_{q_r}k_{q_r} + \frac1{1-\alpha}c_{q_r} > 0$ and define $\bar{k}_{q_r} = c_{q_r} k_{q_r}$. Then, \\
$$2 \lambda_{\min}(W_{q_r}) = 2(1 - \alpha)\bar{k}_{q_r} - c_{q_r}\left(\frac{2\alpha}{1-\alpha}C_{q_r}^2 + C_{q_r} + (1+\alpha)\bar{k}_{q_r} + \frac{1+\alpha}{1-\alpha}c_{q_r} \right),$$
from which it is clear that $\lambda_{\min}(W_{q_r})$ can be made arbitrarily large (and positive) by choosing $c_{q_r}$ and $\bar{k}_{q_r}$ appropriately. Another application of Lemma \ref{posdeflemma} then shows that $W_{q_r} - \frac14 \bar{W}_{q_r, q_j} W_{q_j}^{-1} \bar{W}_{q_r, q_j} \succ 0$, and its minimum eigenvalue can be made arbitrarily large with appropriate choices of $k_{q_j}, k_{\omega_j}, c_{q_j}$ for $j \in \mathcal{I}$.

Now we look at condition 3. First, choose $k_{x_r}, k_{v_r}, c_{x_r}$ such that $W_{x_r} \succ 0$ (this can always be done by, for instance, choosing $c_{x_r}$ sufficiently small). We now wish to show that the remaining subtractive terms can be shrunk arbitrarily. Observe that $\bar{W}_{x_r, q_j}$ and $W_{q_j}^{-1}$ are independent, so that we may force the maximum eigenvalue of $\bar{W}_{x_r, q_j} W_{q_j}^{-1} \bar{W}_{x_r, q_j} \succeq 0$ to be arbitrarily small after (potentially) further shrinking the maximum eigenvalue of $W_{q_j}^{-1}$. Observe that we may write the third term, $\frac14 \bar{W}_{x_r, q_r} (W_{q_r} - \frac14 \bar{W}_{q_r, q_j} W_{q_j}^{-1} \bar{W}_{q_r, q_j})^{-1} \bar{W}_{x_r, q_r}$, in the form $\alpha^2 M^T A M$, where $M$ is independent of $\alpha$ and the terms of $A$ are at most of $O(\frac1{\alpha})$. Hence, we may shrink this term arbitrarily by shrinking $\alpha.$ The fourth and fifth terms are transposes of each other and therefore may be handled simultaneously. Note that find that the maximum eigenvalue is bounded above by
\begin{align*}
    \lambda_{\max}((\bar{W}_{q_r} -& \frac14 \bar{W}_{q_r, q_j} \bar{W}_{q_j}^{-1}\bar{W}_{q_r, q_j})^{-1})||\bar{W}_{q_r, q_j} \bar{W}_{q_j}^{-1} \bar{W}_{x_r, q_j}|| \ ||\bar{W}_{x_r, q_r}|| \\
    \le &\sqrt{\lambda_{\max}(\bar{W}_{q_j}^{-1})} \lambda_{\max}((\bar{W}_{q_r} - \frac14 \bar{W}_{q_r, q_j} \bar{W}_{q_j}^{-1}\bar{W}_{q_r, q_j})^{-1}) ||\bar{W}_{q_r, q_j}|| \ || \bar{W}_{x_r, q_j}|| \ ||\bar{W}_{x_r, q_r}||.
\end{align*}

This term therefore can be arbitrarily shrunk by shrinking the maximum eigenvalue of $\bar{W}_{q_j}$. Moreover, the presence of the norm $||W_{x_r, q_r}||$ also gives us control of the size of the term via $\alpha.$ The final term is handled similarly—we find that the maximum eigenvalue is bounded  by \begin{align*}
   \lambda_{\max}((\bar{W}_{q_r} - &\frac14 \bar{W}_{q_r, q_j} \bar{W}_{q_j}^{-1}\bar{W}_{q_r, q_j})^{-1})||\bar{W}_{q_r, q_j} \bar{W}_{q_j}^{-1} \bar{W}_{x_r, q_j}||^2 \\
    \le &\lambda_{\max}(\bar{W}_{q_j}^{-1}) \lambda_{\max}((\bar{W}_{q_r} - \frac14 \bar{W}_{q_r, q_j} \bar{W}_{q_j}^{-1}\bar{W}_{q_r, q_j})^{-1}) ||\bar{W}_{q_r, q_j}||^2 \ || \bar{W}_{x_r, q_j}||^2,
\end{align*} which again may be shrunk arbitrarily by shrinking the maximum eigenvalue of $\bar{W}_{q_j}^{-1}$. In summation, for sufficiently small $\alpha$ and $\lambda_{\max}(W_{q_j})$ sufficiently large, condition 3 is satisfied and $\mathcal{\bar{W}} \succ 0.$ 

Now, letting $x \in S^6$ such that $x = (x_1, x_2) \in \R^4 \times \R^2$, we get that\\ $\lambda_{\min}(\bar{\mathcal{W}}) \ge \min\{ \lambda_{\min}(Q), \ \lambda_{\min}(P - SQ^{-1}S^T)\} ||x + SQ^{-1}x_2||^2$. Note that $||SQ^{-1}x_2||$ can be made arbitrarily small by increasing the maximum eigenvalue of $Q$, so that $||x + SQ^{-1}x_2||$ can be made arbitrarily close to $1$. Further decomposing $P - SQ^{-1}S^T$, we find that the minimum eigenvalue of $\bar{\mathcal{W}}$ is bounded below by a quantity that can be made arbitrarily close to the minimum of the minimum eigenvalue of the matrices in conditions (1), (2), and (3) above—all of which can be made arbitrarily large. Hence, the minimum eigenvalue of $\bar{\mathcal{W}}$ can be made arbitrarily large.
\hfill$\square$

\begin{remark}
Observe that, in the proof of Theorem \ref{gainsth}, it is not important that the minimum eigenvalue can be made arbitrarily large—we need only have that it is positive. However, in Theorem \ref{gainsthdist} of Section \ref{sec_dist}, we will introduce unstructured bounded disturbances to the problem. In such a case, this fact will be crucial. 
\end{remark}

\subsection{Control design for the unreduced model}

Note that for the design of the geometric controllers $u_j$ we assumed that each quadrotor can generates a thrust along any direction. However, the dynamics of each quadrotor is underactuated since the direction of the total thrust is always parallel to its third body-fixed axis (see Fig \ref{figuav}), despite the magnitude of the total thrust can be arbitrarily changed (recall that the total thrust is given by $u_j = f_j R_j e_3$, being $f_j$ the total thrust magnitude and $R_ie_3$ the direction of the third body-fixed axis). The attitude of each quadrotor is controlled such that the third body-fixed axis becomes parallel to the direction of the control force $u_j$ designed by  \eqref{uparallel} and \eqref{ujperp}.

 The desired direction of the third body-fixed axis for each quadrotor, denoted by $b^{3}_j\in S^2$ is given by $b_{j}^{3} = \frac{u_j}{\|u_j\|}.$ Usually, such expression is considered as a constraint on the desired attitude of each quadrotor. Therefore to solve the dimensionality problem arising after introduce the constraint, and to be able to solve the system, the desired direction of the first body-fixed axis $b^{1}_j(t)\in S^2$ is introduced as a smooth function of time. Since the first body-fixed axis is normal to the third body-fixed axis, one can not reach the arbitrary body axis $b^{1}_{j}$ exactly. The usual strategy in this situation is to project into the plane normal to $b^{3}_j$, and the desired direction of the second body-fixed axis is chosen to obtain an orthonormal frame. That is, the desired attitude for each quadrotor are given by $\displaystyle{
R_{j,d} = \begin{bmatrix}\frac{(\hat b^{3}_j)^2 b^{1}_j}{\|(\hat b^{3}_j)^2 b^{1}_j\|}, &
-\frac{\hat b^{3}_jb^{1}_j}{\|\hat b^{3}_jb^{1}_j\|},& -b^{3}_j\end{bmatrix}}\in SO(3)$.

Using \eqref{Rkinred} , the desired angular velocity for each quadrotor is  $\Omega_{j,d} = (R_{j,d}^T\dot R_{j,d})^\vee\in\mathbb{R}^3$, where ${(\cdot)}^\vee:\mathfrak{so}(3)\to\mathbb{R}^3$ denotes the inverse of the hat map.  Define the tracking error vectors for the attitude of each quadrotors $e_{R_j} = \frac{1}{2}(R_{j,d}^T R_j -R_j^T R_{j,d})^\vee$ and the tracking error vectors for the angular velocity of each quadrotor as $e_{\Omega_j} = \Omega_j - R_j^T R_{j,d}\Omega_{j,d}$.
The thrust magnitude is chosen as the length of $u_j$, projected on to $R_je_3$, and the control moment is chosen as a tracking controller on $SO(3)$, that is, 
\begin{align}
f_j & = u_j\cdot R_j e_3,\label{eqn:fi}\\
M_j & = -\frac{k_R}{\epsilon^2} e_{R_j} -\frac{k_\Omega}{\epsilon} e_{\Omega_j} + \Omega_j\times J_j\Omega_j- J_j (\hat\Omega_j R_j^T R_{j,d}\Omega_{j,d}-R_j^T R_{j,d}\dot\Omega_{j,d}),\qquad j=1,2,\label{eqn:Mi}
\end{align} where $\epsilon,k_R,k_\Omega$ are positive constants. 

Stability of the corresponding controlled systems for the unreduced model can be studied by using singular perturbation theory for the attitude dynamics of quadrotors as in~\cite{SreLeePICDC13,LeeSrePICDC13}. In particular, as a direct application of \cite{Lee-CDC,Lee-TCST}, in the context of Theorem \ref{gainsth} for our particular cooperative transportation task lead to the following result.

\begin{corollary}\label{corol}
Consider the control system defined
by \eqref{xred}-\eqref{Rdynred} and the control inputs designed by  \eqref{eqn:fi} and \eqref{eqn:Mi}. Then, there
exists $\delta > 0$, such that for all $\epsilon< \delta$, the zero equilibrium of the tracking errors\\$(e_{x_r}, e_{v_r}, e_{q_r}, e_{\omega_r}, e_{q_1}, e_{\omega_1}, e_{q_2}, e_{\omega_2}, e_{R_1},e_{\Omega_1},e_{R_2},e_{\Omega_2}) $ is exponentially stable.
\end{corollary}

Now that we have established the exponential tracking of the full reduced model, we wish to connect this back to our original model with elastic cables. This can be done by showing that our system is under the conditions of Theorem 11.2 in $\cite{Khalil}$. 

Before stating the Proposition formally, we introduce some notation and definitions that will make the statement more compact. 

\begin{definition}
The \textit{boundary layer system} for the singular perturbation problem given by \eqref{eqslow}-\eqref{eqfast} is defined as:
$$\frac{\partial r}{\partial \tau} = g(t,x,r+h(t,x), 0),$$where $r := z - h(t,x)$ with $h(t,x)$ as defined by \eqref{heq1}-\eqref{heq2} and $\tau := \frac{t - t_0}{\epsilon}$ for $t_0$ the value of time from which we obtain our initial data. 
\end{definition}

The following Corollary for the exponential stability of the boundary layer system for \eqref{eqslow}-\eqref{eqfast} follows from the case of a single quadrotor transporting a point mass load with an elastic cable (see Lemma 2 in \cite{Elastic}). 

\begin{corollary}\label{boundarylayer}
The boundary layer system for \eqref{eqslow}-\eqref{eqfast} with control inputs $u_j$ and $M_j$ as defined above has an exponentially stable equilibrium point at the origin.
\end{corollary}

Theorem 11.2 in $\cite{Khalil}$ tells us that the trajectories of the original model lie in a neighborhood of the trajectories of the reduced model when the origin of the boundary layer system and the error dynamics of the reduced model are exponentially stable -- which follows immediately from Corollary \ref{boundarylayer} and Theorem \ref{gainsth} above. Formally stated, we have the following Proposition:

\begin{proposition}\label{singular}
Let the control inputs $u_j$ and $M_j$ be defined as above. Denote by $x(t)$ a trajectory of the reduced model \eqref{xred}-\eqref{Rdynred} which converges exponentially to the desired trajectory. Denote by $r(t)$ a trajectory of the boundary layer system which converges exponentially to the origin. Then, there exists a positive constant $\epsilon^\ast$ such that for all $t \ge t_0$ and $0 < \epsilon < \epsilon^\ast$, there exists a unique solution $x(t, \epsilon)$, $z(t, \epsilon)$ of the singular perturbation problem \eqref{eqslow}-\eqref{eqfast} on $[t_0, \infty)$ satisfying 
\begin{align*}
x(t, \epsilon) - x(t) &= \mathcal{O}(\epsilon),\\
z(t, \epsilon) - h(t, x(t)) - r\left(\frac{t - t_0}{\epsilon}\right) &= \mathcal{O}(\epsilon),
\end{align*}
uniformly on $t \in [t_0, \infty)$. Moreover, for $t_1 > t_0$, we have \begin{align*} z(t,\epsilon) - h(t,x(t)) = O(\epsilon),\end{align*} uniformly on $(t_1, \infty)$ for $\epsilon < \epsilon^{\ast\ast} < \epsilon^\ast$.
\end{proposition}

\end{section}

\subsection{Simulation Results}

We now test the results of Theorem \ref{gainsth} with numerical simulations. Simulations of the dynamics for the reduced model were conducted using the proposed controller.

We implement an Euler method for the numerical integration of equations \eqref{xkin2}-\eqref{qjdyn2} with time step $h = 0.002$ and $T=20$ sec. The system parameters were chosen as $m_Q = 0.755, \ m_r = 0.5, \ L_c = L_r = 1.0,$ and $J_Q = \begin{bmatrix} 0.082 & 0 & 0 \\ 0 & 0.0845 & 0 \\ 0 & 0 & 0.1377 \end{bmatrix}$. The gains we consider are $k_{x_r} = 9, \ k_{v_r} = 6, \ k_{q_r} = 2, \ k_{w_r} = 2 \sqrt{2}, \ k_{q_1} = 36, \ k_{w_1} = 12, \ k_{q_2} = 36, \ k_{w_2} = 12$. It should be noted that both smaller and larger gain-sets were found that yielded convergence. This particular set was chosen because it offered an appropriate balance between performance and feasibility, in particular with regards to the convergence of the load position and attitude to their desired trajectories, which were chosen as the \textit{Lissajous curve} described by
$\tilde{x}(t)= (1.2 \sin(0.4 \pi t), \ 4.2 \cos(0.2 \pi t), \ -0.5)$, $\tilde{q}_r(t) = (0, 1, 0)$. The initial conditions were chosen randomly as $x_r(0) = (1, 4.9, -1), \ v_r(0) = (1.2, 0.55, 0.15), \ q_r = (0.24, 0.97, -0.1), \ \omega_r = (0.1, -0.1, 0), \ q_1 = (0.53, 0.63, -0.56), \ \omega_1 = (0, 0, 0), \ q_2 = (0.48, 0.67, -0.56), \omega_2 = (0, 0, 0).$

\begin{figure}[h!] \includegraphics[width=7.0cm]{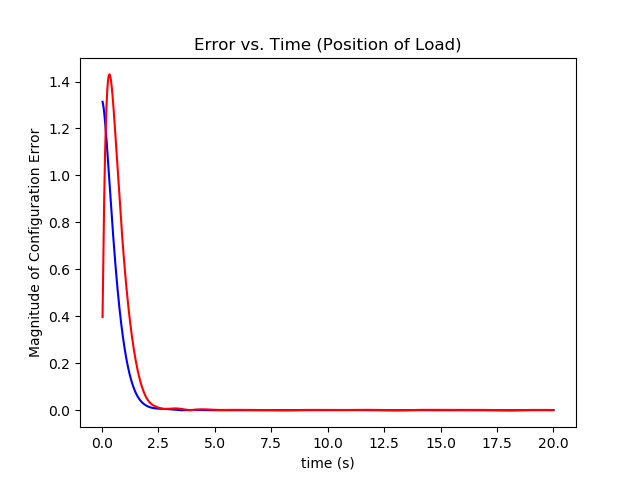}
\includegraphics[width=7.0cm]{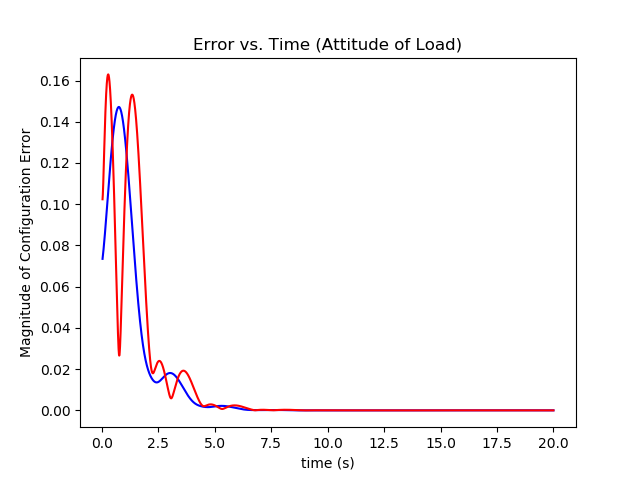}
\includegraphics[width=7.0cm]{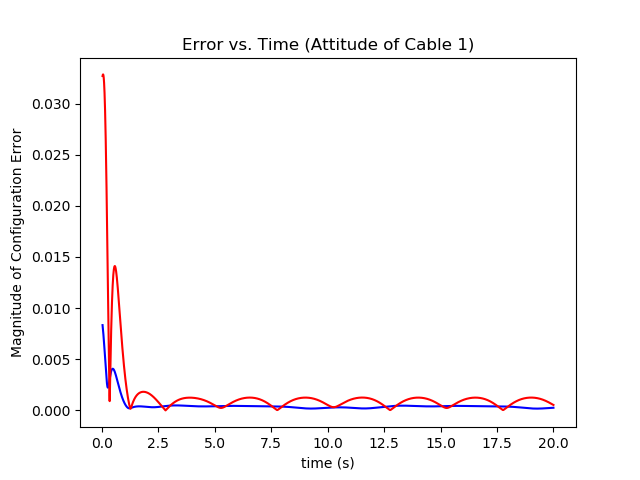}
\includegraphics[width=7.0cm]{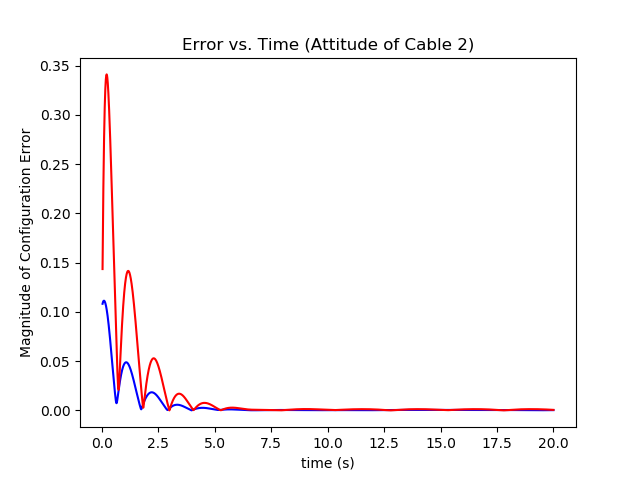}
\includegraphics[width=7.0cm]{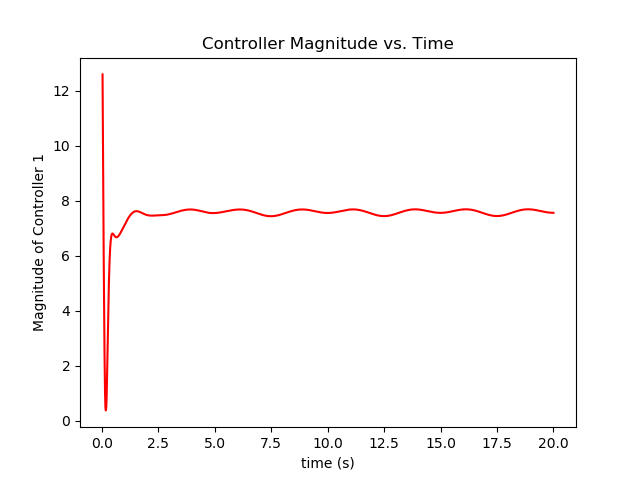}
\qquad \qquad\qquad\quad\includegraphics[width=7.0cm]{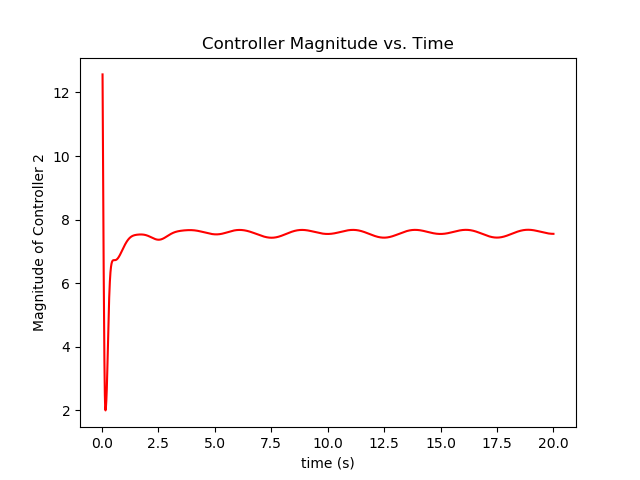}
 \caption{Configuration errors vs. time and control magnitude.}
    \label{fig:plot3}
\end{figure}

Figure \ref{fig:plot3} shows the magnitude of the error functions versus times. In each case, the blue curve represents the error in the position/attitude, and red curve represents the error in the velocity/angular velocity. Additionally, we plot the magnitude of the controllers $||u_1||$ and $||u_2||$ versus time.

\begin{figure}[h!]
    \centering
    \includegraphics[width=7.5cm]{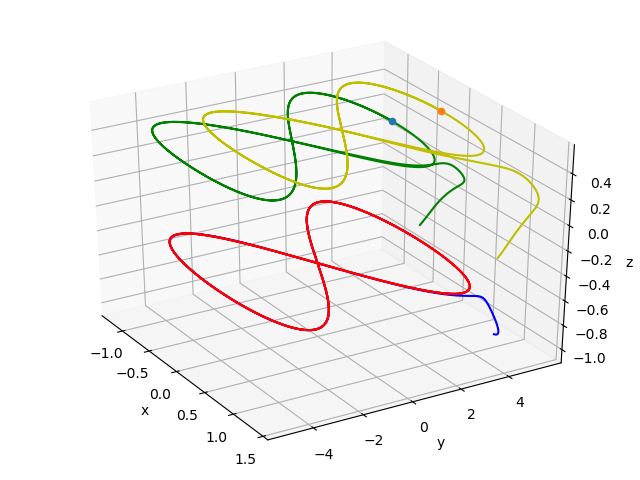}
    \includegraphics[width=7.5cm]{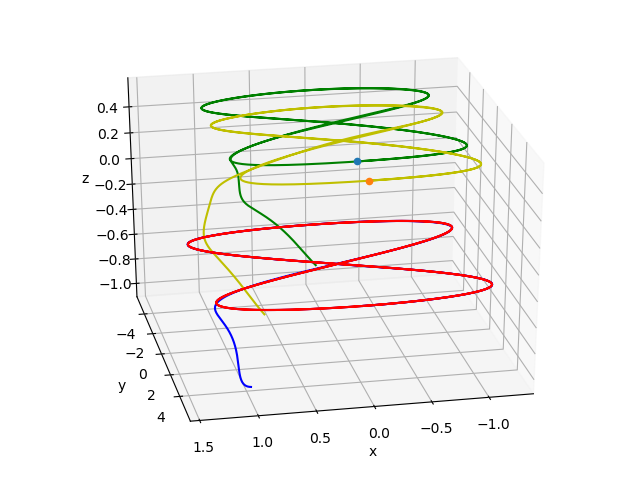}
    \caption{Position curves for the load (blue) and two quadrotors (green and yellow), together with the desired trajectory of the load (red).}
    \label{fig:plot2}
\end{figure}

Note that, as expected, the error asymptotically converge to zero (or a small neighborhood of zero) in all cases, while the controls approach steady state solutions. For illustrative purposes, we also include Figure \ref{fig:plot2}, which plots the load position (blue curve), the desired trajectory for the load (red curve), and the positions of the quadrotors (green and yellow curves) in $3D$ space. The highlighted points along the yellow and green curves mark the final positions of the two quadrotors. 

\begin{section}{Extension of the control design in the presence of (unstructured) disturbances}\label{sec_dist}

While the control scheme developed in Section \ref{sec5} is theoretically sound, in practice there may be practical constraints such as disturbances and measurement errors which cause its failure. This in part is handled by the well-known fact that exponential stability is robust to small disturbances \cite{khalil2}. However, it would additionally be beneficial to understand how the proposed controller behaves when it is subject to larger disturbances. Towards that end, we introduce bounded unstructured disturbances to the reduced dynamical model \eqref{xkin2}-\eqref{qjdyn2}. That is, for some real numbers $\delta x_r, \delta q_j$ for $j \in \mathcal{I}$, we consider the dynamical system
\begin{align}
\dot{x}_r &= v_r, \label{xkin2dis}\\
m_r ( \dot{v}_r + ge_3) &= \mu_1 + \mu_2 + \Delta x_r, \label{xdyn2dis}\\
\dot{q}_r &= \omega_r \times q_r, \label{qrkin2dis}\\
\frac23 m_r L_r \dot{\omega}_r &= q_r \times ( \mu_2 - \mu_1 ) + \Delta q_r, \label{qrdyn2dis} \\
\dot{q}_j &= \omega_j \times q_j, \label{qjkin2dis} \\
m_Q L_c \dot{\omega}_j &= m_Q q_j \times ( \ddot{x}_r + (-1)^j L_r \ddot{q}_r + ge_3) - q_j \times u_j + \Delta q_j, \label{qjdyn2dis}
\end{align} where for $j \in \mathcal{I}$, $\Delta x_r, \Delta q_r, \Delta q_j$ are unstructured disturbances satisfying $||\Delta x_r|| \le \delta x_r$ and $||\Delta q_j|| \le \delta q_j$. Defining the controllers, configuration errors, and Lyapunov candidate $V$ as in section \eqref{sec5}, the bounds on the Lyapunov function remain the same. That is, we have $\frac12 z_{q_j}^T \underbar{P}_{q_j} z_{q_j} \le V_{q_j} \le \frac12 z_{q_j}^T \bar{P}_{q_j} z_{q_j}.$ However, some additional terms appear in its derivative (last three terms in the following inequality):

\begin{align*}
\dot{V}  \le & -(k_{v_r} - c_{x_r}) ||e_{v_r}||^2 + c_{x_r} k_{v_r} ||e_{v_r}|| ||e_{x_r}|| - c_{x_r}k_{x_r}||e_{x_r}||^2 + \frac1{m_r}(||e_{v_r}|| + c_{x_r} ||e_{x_r}||) ||Y|| \\
&  -(k_{\omega_r} - c_{q_r})||e_{\omega_r}||^2 + c_{q_r}(k_{\omega_r} + C_{q_r})||e_{\omega_r}||||e_{q_r}|| - c_{q_r} k_{q_r} ||e_{q_r}||^2 - c_{q_j} k_{q_j} ||e_{q_j}||^2\\ &+ \frac2{3 m_r L_r} (||e_{\omega_r}|| + c_{q_r} ||e_{q_r}||) ||Y||   -(k_{\omega_j} - c_{q_j}) ||e_{\omega_j}||^2 + c_{q_j} (k_{\omega_j} + C_{q_j} ) ||e_{q_j}|| ||e_{\omega_j}|| \\
& + E_{x_r}^T (e_{v_r} + c_{x_r} e_{x_r}) + E_{q_r}^T (e_{\omega_r} + c_{q_r} e_{q_r}) + E_{q_j}^T (e_{\omega_j} + c_{q_j} e_{q_j}),
\end{align*} where $E_{x_r}, E_{q_r},$ and $E_{q_j}$ are defined as $E_{x_r} = \frac{\Delta x_r}{m_r}, 
E_{q_r} = \frac{3 \Delta q_r}{2m_r L_r}, 
E_{q_j} = \frac{\Delta q_j}{m_Q L_c}$. We may then write this is the form $\dot{V} \le -z^T \mathcal{W} z + E^T z,$ where 
$E := \begin{bmatrix}
\frac{\delta x_r}{m_r}, \frac{\delta x_r}{m_r}, \frac{3\delta q_r}{2m_r L_r}, \frac{3\delta q_r}{2m_r L_r}, \frac{\delta q_j}{m_Q L_c}, \frac{\delta q_j}{m_Q L_c} 
\end{bmatrix}$ and $\mathcal{W}$ is as in section \ref{sec5}.
From Young's Inequality, we have $E^T z \le \frac{||E||^2}{4 \epsilon} + \epsilon ||z||^2$ for any $\epsilon > 0$. Hence,
$$\dot{V} \le -z^T (\mathcal{W} - \epsilon I) z + \frac{||E||^2}{4 \epsilon}.$$

Fix $\epsilon > 0.$ Analogous to before, we replace the matrix $\mathcal{W} - \epsilon I$ with its symmetric part $\bar{\mathcal{W}}^\ast := \frac12(\mathcal{W} + \mathcal{W}^T) - \epsilon I$. Observe that:

\begin{align*}
    \lambda_{\min}(\mathcal{\bar{W}}^\ast) &= \min_{x \in S^6} \{ x^T (\mathcal{\bar{W}} - \epsilon I) x \} = \min_{x \in S^6} \{ x^T \mathcal{\bar{W}} x - \epsilon ||x||^2 \} = \min_{x \in S^6} \{ x^T \mathcal{\bar{W}} x \} - \epsilon \\
    &= \lambda_{\min}(\mathcal{\bar{W}}) - \epsilon 
\end{align*}

From the proof of Theorem \ref{gainsth}, $\lambda_{\min}(\mathcal{\bar{W}})$ can be made arbitrarily large by choosing gains appropriately. Hence, we can choose them so that $\bar{\mathcal{W}}^\ast \succ 0$, with an arbitrarily large minimum eigenvalue. We then have the following inequalities:

\begin{align*}
\lambda_{\min}(\underbar{P})||z||^2 \le V \le \lambda_{\max}(\bar{P})||z||^2, \\
\dot{V} \le -\lambda_{\min}(\bar{\mathcal{W}}^\ast)||z||^2 + \frac{||E||^2}{4\epsilon}.
\end{align*}

This implies that $\dot{V} \le -\frac{\lambda_{\min}(\bar{\mathcal{W}}^\ast)}{\lambda_{\max}(\bar{P})}V + \frac{||E||^2}{4\epsilon}$, so that $\dot{V} < 0$ when $V > \frac{\lambda_{\max}(\bar{P})}{\lambda_{\min}(\bar{\mathcal{W}}^\ast)}\frac{||E||^2}{4\epsilon}:= d_1 > 0.$ Clearly, $d_1$ can be arbitrarily small by making $\lambda_{\min}(\bar{\mathcal{W}}^\ast)$ sufficiently large. If we now define the set $S_{r} := \{z \in \mathcal{D} : V(z) < r \}$, where $r$ is some real number, then any trajectory starting in the open set $\mathcal{D} \setminus \bar{S}_{d_1}$ will converge exponentially to the region $\bar{S}_{d_1}$, where $\bar{S}_{d_1}$ denotes the topological closure of $S_{d_1}$. Since $V$ is continuous and positive, $\bar{S}_{d_1}$ is some closed neighborhood of the origin that can be made arbitrarily small (by shrinking $d_1$). We formalize this result with the following theorem:

\begin{theorem}\label{gainsthdist}
Consider the control system with disturbances defined by equations \ref{xkin2dis} - \ref{qjdyn2dis} with control inputs \ref{mutilde}, \ref{mutildeq}, and \ref{ujperp}. For sufficiently small $\alpha$, there exists control gains $k_{x_r}$, $k_{v_r}$, $k_{q_r}$, $k_{\omega_r}$, $k_{q_1}$, $k_{\omega_1}$, $k_{q_2}$, and $k_{\omega_2}$ such that the zero equilibrium of the tracking errors $e_{x_r}$, $e_{v_r}$, $e_{q_r}$, $e_{\omega_r}$, $e_{q_1}$, $e_{\omega_1}$, $e_{q_2}$ and $e_{\omega_2}$ are uniformly ultimately bounded. Moreover, the ultimate bound can be made arbitrarily small. 
\end{theorem}

\begin{remark}
We will not extend this result to include the quadrotor attitude kinematics and dynamics here, nor relate it back to the original elastic model. Such a task follows essentially the same strategy as it did in Section \ref{sec5}, just replacing exponential stability with uniform ultimate boundedness where appropriate.  
\end{remark}

\end{section}

\begin{section}{Conclusion and future work}
In this paper, we propose a model and geometric trajectory tracking controller for the cooperative task of two quadrotor UAVs transporting a rigid bar via inflexible elastic cables. This is handled in three stages: (i) Reduction of the model to that of a similar model with inelastic cables. We accomplish this by assuming sufficient stiffness and damping of the cables (a realistic condition in applications) and utilizing the results of singular perturbation theory. (ii) Design of a geometric tracking controller in the reduced model. Lyapunov analysis is used to find sufficient conditions for stability, and Theorem \ref{gainsth} proves the existence of gains satisfying these conditions for sufficiently small initial errors in the cable attitudes. (ii) We show that—under the same control law—trajectories of the original (elastic) model converge uniformly to the trajectories of the reduced model as the stiffness and damping of the cables approach infinity. We also extended the proposed approach to design a control law for the case when the system under study is subject to unstructured bounded disturbances.

We are currently working to add uncertainties in order to further explore the robustness of the proposed controller. We also plan to study the construction of force variational integrators in optimal control problems, in a similar fashion to \cite{leo1} and \cite{leo2}, dynamic interpolation problems \cite{interpolation}, and obstacle avoidance problems \cite{variational} for the cooperative mission. Finally, note also that in our model, cables are attached to the center of each quadrotor. It would be interesting to shift those attachment points and see how to deal with the resulting coupled systems—instead of a decoupled system as it is in this paper—as well as heterogeneous agents.

\end{section}

\section*{Acknowledgments}
Leonardo Colombo (leo.colombo@icmat.es) and Jacob Goodman (jacob.goodman@icmat.es) conduct their research at Instituto de Ciencias Matematicas
(CSIC-UAM-UC3M-UCM), Calle Nicolas Cabrera 13-15, 28049, Madrid, Spain. The project that gave rise to these results received the support of a fellowship from ”la Caixa” Foundation (ID 100010434). The fellowship codes are LCF/BQ/PI19/11690016 and LCF/BQ/DI19/11730028. The authors also acknowledge financial support from the Spanish Ministry of Science and Innovation, under grants PID2019-106715GB-C21, MTM2016-76702-P, ``Severo Ochoa Program for Centres
of Excellence in R\&D'' (SEV-2015-0554) and from the Spanish National
Research Council, through the ``Ayuda extraordinaria a Centros de Excelencia Severo Ochoa''(20205CEX001). All the results of the paper are original and have not been presented nor submitted to a conference. The authors also acknowledge Jos\'e Angel Acosta from University of Seville, Spain, for fruitful discussions on the inclusion of disturbances into our model. The authors also wish to thank Ravi Banavar from IIT-Bombay and D.H.S. Maithripala from University of Peradeniya, for helpful discussions about the modelling of cooperative tasks with quadrotors.


\begin{thebibliography}{20}
\providecommand{\newblock}{\relax}

\bibitem{abloch}A. M. Bloch. Nonholonomic mechanics and control.   Springer-Verlag, Second Edition, 2015.
\bibitem{interpolation}  A. Bloch, M. Camarinha and L. J. Colombo. Dynamic interpolation for obstacle avoidance on Riemannian manifolds. International Journal of Control,
pages 1-22, doi:10.1080/00207179.2019.1603400, https://doi.org/10.
1080/00207179.2019.1603400. Preprint available at. arXiv:1809.03168
[math.OC]
\bibitem{variational} A. Bloch, M. Camarinha, L. Colombo. Variational obstacle avoidance on Riemannian manifolds. in Proceedings of the IEEE International Conference on Decision
and Control, 2017, pp. 146-150.
\bibitem{amit1}
N.A. Chaturvedi, A.K. Sanyal, N.H. McClamroch. Rigid-body attitude control.
IEEE control systems magazine 31 (3), 30-51, 2011.
\bibitem{leo1} L. Colombo, S. Ferraro, D. Martn de Diego. Geometric integrators for higherorder variational systems and their application to optimal control. J. Nonlinear
Sci. 26 (2016), no. 6, 1615-1650.
\bibitem{leo2} L. Colombo, F. Jimenez, and D. De Diego, “Variational integrators ´
for mechanical control systems with symmetries,” Journal of Computational Dynamics, vol. 2, no. 2, pp. 193–225, 2015.
\bibitem{gas} M. Gassner, T. Cieslewski, D. Scaramuzza. Dynamic collaboration without communication: Vision-based cable-suspended load transport with two quadrotors, IEEE International Conference on Robotics and Automation, 5196-5202, 2017.


\bibitem{GF} I. M. Gelfand and S. V. Fomin. Calculus of variations. Revised English edition translated and edited by Richard A. Silverman. Prentice-Hall, Inc., Englewood Cliffs, N.J., 1963, pp. vii+232

\bibitem{god} F. A. Goodarzi and T. Lee, “Dynamics and control of quadrotor
uavs transporting a rigid body connected via flexible cables,” in 2015
American Control Conference, 2015, pp. 4677-4682.

\bibitem{SPT} M. Holmes. Introduction to perturbation methods, volume 20. Springer Science $\&$ Business Media, 2012.
\bibitem{HSS}
D. D. Holm, T. Schmah, C. Stoica,\textit{ Geometric mechanics and symmetry,} Oxford University Press, 2009.

\bibitem{amit3} M. Izadi, A.K. Sanyal, R.R. Warier. Variational Attitude and Pose Estimation Using the Lagrange-d’Alembert Principle.
2018 IEEE Conference on Decision and Control (CDC), 1270-1275.
\bibitem{amit4} M. Izadi, A.K. Sanyal. Rigid body pose estimation based on the Lagrange–d’Alembert principle. Automatica 71, 78-88, 2016.
\bibitem{jo} S. Joshi and C. D. Rahn, “Position control of a flexible cable gantry
crane: theory and experiment,” in American Control Conference,
Proceedings of the 1995, vol. 4, 1995, 2820-2824.
\bibitem{Khalil}
H. K. Khalil, \textit{Nonlinear Systems,} Prentice-Hall, Inc., Third Edition, 2002.
\bibitem{khalil2} H. K. Khalil, \textit{Nonlinear Control,} Pearson Education Limited, 2015.
\bibitem{LLM}  T. Lee, M. Leok, and N. McClamroch, \textit{Geometric tracking control of a
quadrotor UAV on SE(3)}, in Proc. 49th IEEE Conf. Decision Control,
Dec. 2010, pp. 5420-5425.
\bibitem{LeeSrePICDC13}
T.~Lee, K.~Sreenath, and V.~Kumar, Geometric control of cooperating multiple
  quadrotor {UAV}s with a suspended load.  \textit{Proceedings of the IEEE
  Conference on Decision and Control}, 2013, pp. 5510--5515.

\bibitem{Lee-TAC} T. Lee. Global Exponential Attitude Tracking Controls on $SO(3)$. IEEE Transactions on Automatic Control, 60(10):2837-2842, 2015.
\bibitem{Lee-CDC}
T. Lee. Geometric control of multiple quadrotor UAVs transporting a cable-suspended rigid body. in Proc. 53rd IEEE Conf. Decision Control, Dec. 2014, pp. 6155-6160.
\bibitem{Lee-TCST}
T. Lee. Geometric control of quadrotor UAVs transporting a cable-suspended rigid body. IEEE Transactions on Control Systems Technology, vol. 26, no. 1, pp. 255-264, 2018.
\bibitem{MR}
J. E. Marsden, T. S. Ratiu, \textit{Introduction to Mechanics and Symmetry: A Basic Exposition of Classical Mechanical Systems,} Springer-Verlag, New York, 1999.
\bibitem{Elastic}
P. Kotaru, G. Wu and K. Sreenath, "Dynamics and control of a quadrotor with a payload suspended through an elastic cable," 2017 American Control Conference (ACC), 2017, 3906-3913.

\bibitem{MazKonJIRS10}
I.~Maza, K.~Kondak, M.~Bernard, and A.~Ollero. "Multi-UAV cooperation and
  control for load transportation and deployment",. \textit{Journal of Intelligent and Robotic Systems}, vol.~57, pp. 417--449, 2010.


\bibitem{muga} D. H. S. Maithripala, J. Berg, W. Dayawansa, Almost-global tracking of simple mechanical systems on a general class of lie groups,
IEEE Transactions on Automatic Control 51 (2) (2006) 216-225.
\bibitem{muga2} D. H. S. Maithripala, W. P. Dayawansa, J. M. Berg, Intrinsic observer-based stabilization for simple mechanical systems on lie groups,
SIAM Journal on Control and Optimization 44 (5) (2005) 1691–1711. 
\bibitem{MicFinAR11}
N.~Michael, J.~Fink, and V.~Kumar, ``Cooperative manipulation and
  transportation with aerial robots,'' \textit{Autonomous Robots}, vol.~30, pp.
  73--86, 2011.


\bibitem{PalCruIRAM12}
I.~Palunko, P.~Cruz, and R.~Fierro, ``Agile load transportation,'' \textit{IEEE
  Robotics and Automation Magazine}, vol.~19, no.~3, pp. 69--79, 2012.

\bibitem{pedro1}  P. Pereira and D. V. Dimarogonas. Pose and Position Trajectory Tracking for Aerial Transportation of a Rod-Like Object, Automatica, Vol. 109, pp. 108547, 2019.
\bibitem{pedro2} P. Pereira and D. V. Dimarogonas. Pose Stabilization of a Bar Tethered to Two Aerial Vehicles, Automatica, Vol. 112, pp. 108695, 2020.
\bibitem{amit2}A Sanyal, N Nordkvist, M Chyba. An almost global tracking control scheme for maneuverable autonomous vehicles and its discretization. IEEE Transactions on Automatic control 56 (2), 457-462, 2011.
\bibitem{SreLeePICDC13}
K.~Sreenath, T.~Lee, and V.~Kumar, ``Geometric control and differential
  flatness of a quadrotor {UAV} with a cable-suspended load,'' in
  \textit{IEEE Conference on Decision and Control}, 2013,
  2269--2274.
\bibitem{jose}S. Thapa, H. Bai, J. Acosta Cooperative Aerial Manipulation with Decentralized Adaptive Force-Consensus Control. Journal of Intelligent \& Robotic Systems 97 (1), 171-183, 2020.
\bibitem{wu} G. Wu and K. Sreenath, “Geometric control of quadrotors transporting
a rigid-body load,” in IEEE Conference on Decision and Control, Los
Angeles, CA, Dec. 2014, pp. 6141-6148.

\end{thebibliography}
\end{document}